\documentclass[10pt, a4paper, twoside,reqno]{amsart}
\usepackage[dvipsnames]{xcolor}
\usepackage{pgfplots}
\usepackage{siunitx}
\usepackage{algorithm2e}

\newcommand{\diagIndtwo}[3]      {\mathrm{diag}_{#1}^{#2}\!\left(#3\right)}		
\newcommand{\bma}[2]       {\left[\begin{array}{#1}#2\end{array}\right]}

\newcommand{\trans}        {\top}			
\newcommand{\iter}        {\kappa}	
\newcommand{\w}        {q}

\newcommand{\Pocp} {\mathcal{P}_{\mathrm{ocp}}}
\newcommand{\focp} {f_{\mathrm{ocp}}}
\newcommand{\Pocpu} {\mathcal{P}_{\mathrm{ocp,u}}}
\newcommand{\focpu} {f_{\mathrm{ocp,u}}}
\newcommand{\PSQPGN} {\mathcal{P}_{\mathrm{SQP}}^\iter}
\newcommand{\fSQPGN} {f_{\mathrm{ocp}}^{\mathrm{GN}}}

\newcommand{\focpuGN} {f_{\mathrm{ocp,u}}^{\mathrm{GN}}}
\newcommand{\frF} {R}
\newcommand{\uPocpu} {u_{\mathcal{P}_{\mathrm{ocp,u}}}}
\newcommand{\ul}{u_{l_1}}
\newcommand{\uli}{u_{i,l_1}}
\newcommand{\rFl}{r_{F,l_1}}
\newcommand{\rl}{r_{l_1}}
\newcommand{\rsl}{r_{s,l_1}}

\usepackage{verbatim}

\usepackage{fancyhdr}
\usepackage{mathabx}

\makeatletter

\def\section{\@startsection{section}{1}%
	\z@{.7\linespacing\@plus\linespacing}{.5\linespacing}%
	{\normalfont \Large\scshape\centering}}

\def\subsection{\@startsection{subsection}{2}%
	\z@{.5\linespacing\@plus.7\linespacing}{.5\linespacing}%
	{\normalfont\large\bfseries}}

\def\subsubsection{\@startsection{subsubsection}{3}%
	\z@{.5\linespacing\@plus.7\linespacing}{.5\linespacing}%
	{\normalfont\itshape}}

\definecolor{darkblue}{rgb}{0.0, 0.0, 0.45}

\usepackage{hyperref}
\hypersetup{colorlinks	= true,
	raiselinks	= true,
	linkcolor	= darkblue, 
	citecolor	= Mahogany,
	urlcolor	= ForestGreen,
	plainpages	= false}

\allowdisplaybreaks
\date{\today}

\usepackage{graphicx, psfrag}
\usepackage{amsthm}
\usepackage{amsmath} 
\usepackage{amssymb}  
\usepackage{color}
\usepackage{enumerate}
\usepackage{subfigure}
\usepackage{multirow}
\usepackage{booktabs}

\usepackage{tikz}
\usetikzlibrary{calc,patterns,decorations.pathmorphing,decorations.markings}
\usepackage{xspace}

\newtheorem{myrem}{Remark}
\newtheorem{mycor}{Corollary}

\newtheorem{mythe}{Theorem}
\newtheorem{mylem}{Lemma}

\newtheorem{myprop}{Proposition}

\renewcommand{\qedsymbol}{$\blacksquare$}

\newenvironment{IEEEproof}[1][\hspace{2.3mm} \textbf{Proof}]{\itshape \begin{trivlist} 
\item[\hskip \labelsep {\textbf{ #1}:}]}{\qedsymbol \end{trivlist} }

\DeclareMathOperator*{\argmin}{arg\,min}

\title[Control of an Architectural Cable Net Geometry]{Control of an Architectural Cable Net Geometry}        
   \author{Yvonne R.\ St\"urz, Manfred Morari, Roy S.\ Smith }

\thanks{This research was supported by the NCCR Digital Fabrication, funded by the Swiss National Science Foundation (NCCR Digital Fabrication Agreement \# 51NF40-141853). }
\thanks{
Yvonne R.\ St\"urz and Roy S.\ Smith are with the Automatic Control Laboratory, ETH Zurich, Physikstrasse 3, 8092 Zurich, Swizerland. (e-mail: \{stuerzy, rsmith\}@control.ee.ethz.ch). }
\thanks{
Manfred Morari is with the Department of Electrical and Systems Engineering, University of Pennsylvania, Philadelphia 19104, United States. (e-mail: morari@seas.upenn.edu). }

\begin{document} 

\maketitle

\begin{abstract}
Doubly curved thin concrete shells are very efficient building structures, suitable for light-weight construction because of their high structural stability. In the process of constructing such shells, an efficient innovative flexible formwork that is based on a cable network can be used instead of the costly conventional timber formwork. To guarantee the structural properties of such a shell, the desired form, that is designed and optimized in advance, needs to be precisely achieved. The sensitivity of the flexible cable net formwork to fabrication tolerances and uncertainties makes high accuracy challenging. We propose a new construction method where the form of the cable net structure is measured and controlled in a feedback loop during its construction. Two models based on a force and on an energy approach are reviewed and their equivalence is shown. An efficient control algorithm, which is based on a variant of Sequential Quadratic Programming, and which guarantees feasibility at every iteration, is derived. Based on mild assumptions on the cable net, global convergence to a stationary point is shown. For practical applicability, an extension of the control algorithm for computing sparse input vectors is given.  
Experimental results on a cable net formwork prototype for a shell roof structure are presented to demonstrate the control performance. 
\end{abstract}



\section{Introduction}
\label{sec:intro}

Because of their curvature, doubly curved thin concrete shells can be designed with a high stiffness and stability. Loads induce tension and compression forces rather than flexion. Because of their properties, shell elements can span large areas using comparatively little material. Mostly through requiring less concrete, a significant amount of energy can be saved compared to conventional building structures. In addition to their structural advantages, shells are also interesting from an architectural point of view as their doubly curved form enables new aspects of design and expression in buildings, \cite{DEstreeSterk2006}. 

In the construction process of shell structures a so-called formwork is needed as a supporting structure on which to pour or spray the concrete. Conventional formwork which has been used to date is very labor-, material- and time-intensive, as it consists of a large number of tediously manufactured, non-reusable customized timber elements. To overcome this drawback, an innovative flexible formwork can be used which consists of a net of cables or rods and a fabric layer on top, \cite{Veenendaal2014a}\,, \cite{Veenendaal2012}\,, \cite{Torsing2012}\,. The net is pre-stressed such that the weight of the concrete deforms it to the final designed shell form \cite{Mele2010}\,. The tension forces and the weight of both the net and the concrete are supported by a rigid frame at the boundaries where the net is fixed. This new kind of formwork is beneficial in many aspects. Through the standardization and re-use of elements, the amount of material required and the waste are reduced. Furthermore, the construction of the formwork is much faster and therefore less expensive, which could enable the construction of a larger number of doubly curved shell structures in the future. 

The mechanical properties of the concrete shell, such as buckling stability, depend critically on its form, which is the result of an optimization-based design process. 
To satisfy the accuracy requirements of the shell, the cable net tolerances are very tight. 
However, because of uncertainties in the material behavior and fabrication tolerances of the cable net and of the frame, the desired form of the cable net is in general not precisely achieved if typical construction methods are used. This is shown in \cite{Veenendaal2014} in an experiment for a small-scale simple shell prototype. 
Therefore, control of the form of the cable net is introduced, to achieve the designed target form as closely as possible. 

For the innovative cable net formwork application, both a model based on force equilibrium and based on an energy approach are introduced in \cite{Stuerz2016a}. 
Based on these models, a fast algorithm is presented for the control of the form of the cable net during its construction. It is based on a variant of Sequential Quadratic Programming (SQP) with feasible iterates. 
In \cite{Stuerz2016b}, methods for the identification of important parameters of the cable net system are presented. 

In this paper, the following extensions to \cite{Stuerz2016a} and \cite{Stuerz2016b} are made. 
The two model formulations introduced in \cite{Stuerz2016a} are shown to be equivalent. Based on this result, a convergence proof is given for the control algorithm presented in \cite{Stuerz2016a}. 
The algorithm is extended to compute sparse input vectors. 
It enables the practical implementation of feedback corrections of the structure in an acceptable amount of time. This is especially important for large-scale systems without fully automated, or with even manual actuation. 
In contrast to \cite{Stuerz2016a}, where a validation of the control algorithm is given only in terms of simulation results, this paper presents experimental results on a 1:4 scale prototype of the cable net system for the so-called HiLo Roof (High Performance, Low Energy) \cite{HiLo2016}. Details of this prototype are given in \cite{Liew2017}. 

\begin{figure} 
\centering
\includegraphics[width=0.6\columnwidth]{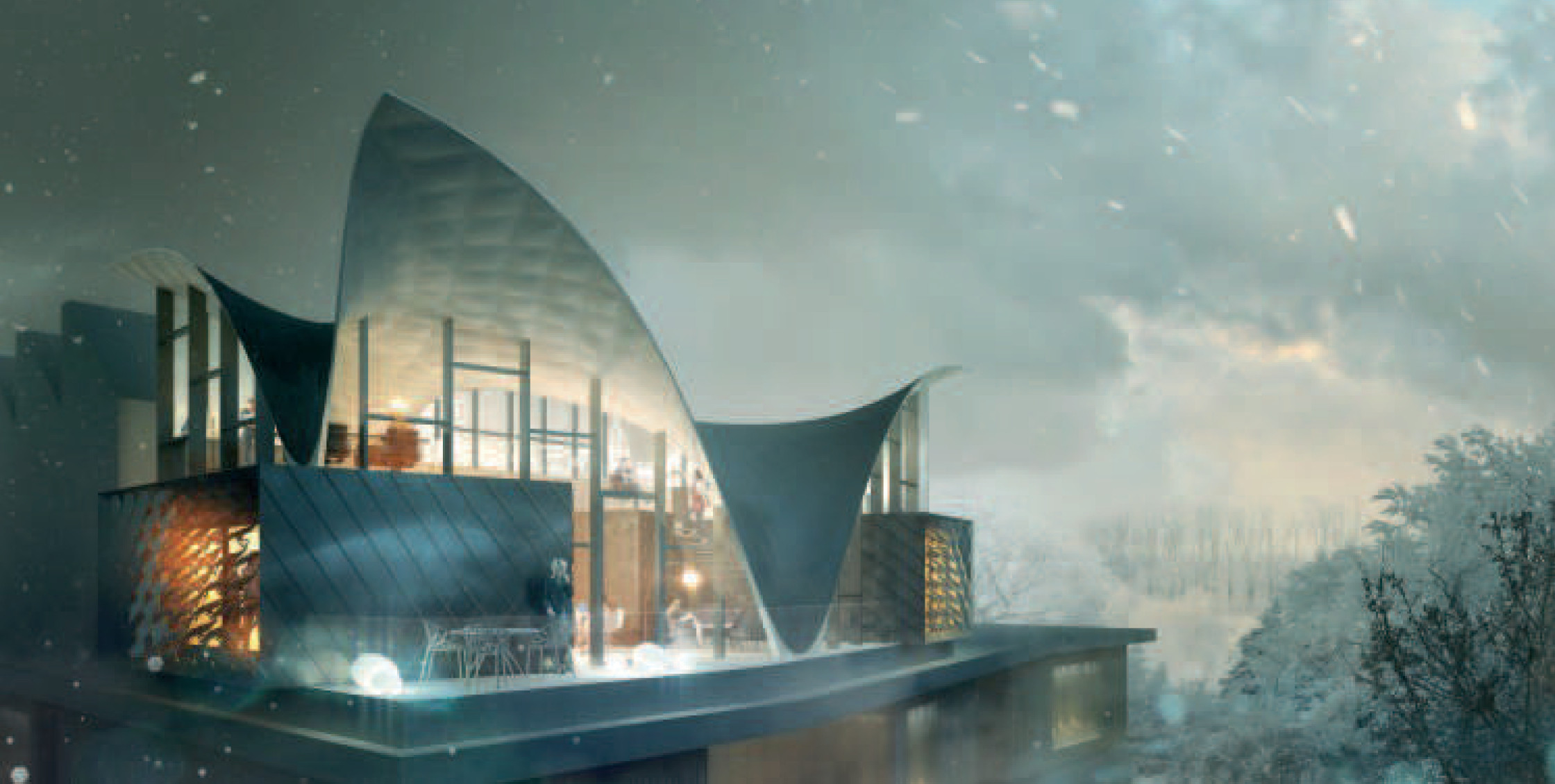}
\caption{Rendering of the HiLo research and innovation unit on the NEST building at Empa, D\"ubendorf, Switzerland. (image by Supermanoevre, Doug and Wolf) \label{fig:HiLo}}
\end{figure}

The HiLo-Roof is part of a research and innovation unit on the so-called NEST-building \cite{Nest2016}\,. The latter is a demonstrator building hosting different research experiments on the campus of the Swiss Federal Laboratories for Materials Science and Technology (Empa), in D\"ubendorf, Switzerland, \cite{Empa2016}\,. The HiLo unit is planned as a duplex penthouse apartment and the roof will have a span of approximately $16 \mathrm{m} \times 9 \mathrm{m}$ and a maximum height of $6.5 \mathrm{m}$. It will be a doubly curved thin concrete shell structure, \cite{Veenendaal2015b}\,, which is depicted in Figure~\ref{fig:HiLo}. For its construction in 2019, a cable net and fabric formwork is going to be used, and the cable net component is planned to be controlled in its form on the construction site. 

The paper is structured as follows. 
Section~\ref{sec:CableNet} gives a description of the cable net formwork, followed by the mathematical model formulations of its equilibrium states in 
Section~\ref{sec:model}. The algorithm for controlling the form is given in 
Section~\ref{sec:Control}. A convergence proof of the control algorithm is given in 
Section~\ref{sec:convergence_GN}. 
Section~\ref{sec:Experiments} gives experimental results conducted on the 1:4 cable net prototype of the HiLo Roof. 

\subsection{Notation} 
We denote a block-diagonal matrix $D$ of submatrices $D_1,...,D_N$ by $D = \diagIndtwo{i=1}{N}{D_i}$. 
The $n \times m$ matrix of all zeros is denoted as $0_{n \times m}$. If clear from the context, the indices are dropped. 
For a function $f(x): \mathbb{R}^n \mapsto \mathbb{R}$, we denote the gradient by $\nabla_{x} f = [\frac{\partial f}{\partial x_1} ~ \cdots ~ \frac{\partial f}{\partial x_n}]^\trans$. For a function $h(r,u): \mathbb{R}^{n+p} \mapsto \mathbb{R}^m$, the Jacobian is denoted by 
\begin{equation}
\begin{aligned}
\nabla_{(r,u)} h &= \bma{@{}c@{\,\,\,} | @{\,\,\,}c@{}}{ 
  \nabla_{r} h & \nabla_{u} h  }  
  = \bma
{@{}c@{\,\,}c@{\,\,}c@{\,\,\,} | c@{\,\,}c@{\,\,}c@{}}{ 
 \frac{\partial h_1}{\partial r_1} & \cdots & \frac{\partial h_1}{\partial r_{n}}  & \frac{\partial h_1}{\partial u_1}  & \cdots & \frac{\partial h_1}{\partial u_{p}}  \\
\vdots &  & \vdots  & \vdots &  & \vdots \\
 \frac{\partial h_{m}}{\partial r_1} & \cdots & \frac{\partial h_{m}}{\partial r_{n}} & \frac{\partial h_{m}}{\partial u_1}  & \cdots & \frac{\partial h_{m}}{\partial u_{p}}  
},
\end{aligned}
\end{equation}
where $\nabla_{({r},{u})}$ denotes the partial derivatives with respect to ${r}$ and ${u}$ and $\nabla_{r} h$ and $\nabla_{u} h$ are referred to as partial Jacobians. 
We use 
$\nabla_{({r},{u})}  {h}({r}^\iter \,,{u}^\iter)$ to denote the Jacobian evaluated at the point $({r}^\iter,{u}^\iter)$.  
The weighted $L_2$-norm with $Q$ being a weighting matrix is denoted as $\left\Vert x \right\Vert_{Q}^2 = x^\trans Q \, x$.

\section{Cable Net Description}
\label{sec:CableNet}
The cable net is a pre-stressed, pin-jointed structure, which can be seen as a special form of tensegrity structure 
\cite{Skelton}, \cite{Tur2008}, which consists of only cables and has no struts or bars. This kind of tensegrity structure was introduced as spider-web in \cite{Connelly1996}. 
Tensegrity structures have been an active field of research considering for example applications of deployable structures or lightweight structures in different fields such as aerospace or biology. Control of these structures has been proposed for example with the goal of active vibration reduction, where forces or displacements, are applied to actuate the dynamics of the structures, \cite{Chan2004}. 
In \cite{VandeWijdeven2005}, \cite{Wroldsen2009}, \cite{Skelton}, 
the equilibria of the structures are changed by changing the physical parameters of the structure. 
In \cite{Wroldsen2009}, \cite{Skelton}, a Lyapunov-based control for the dynamics of a specific class of tensegrity systems, consisting of a rigid bar connected to strings, is proposed. By changing the initial lengths of the strings as inputs, the rigid-body dynamics of the bar is controlled.  

The control methodology in this work is to guide the cable net system to a different equilibrium which minimizes the error between its actual and desired form. 
This is done by changing the boundary lengths of the net. 
In the control algorithm, we exploit properties of the mathematical model of the cable net equilibria in order to efficiently compute the control inputs. 

This section reviews the cable net, its construction steps and its function as a formwork. The mathematical description in terms of variables, fixed parameters and control inputs is introduced based on a graph-theoretical approach, \cite{Roth1995}, \cite{Stuerz2016a}.

\subsection{Physical Description of the Cable Net}
The innovative flexible formwork is a network of cable elements onto which a fabric membrane is placed. The boundary edges of the net are attached to a rigid frame. 

The design process of the shell comprises several steps. 
A detailed presentation can be found in \cite{Mele2010}\,, \cite{Veenendaal2014}\,, \cite{Veenendaal2014a}\,, \cite{Veenendaal2015}\,. First, the form of the shell is designed, taking into account multiple criteria such as boundary conditions, head clearance, esthetic and design aspects, buckling stability and other mechanical properties. Then the cable net topology is chosen and mapped onto it. Via a best-fit-optimization, the desired force distribution of the cable net loaded by the concrete is obtained. From this final tensioned and loaded state of the cable net, the initial state of the unloaded pre-stressed cable net, i.e. without the concrete, is obtained in terms of its form and its tension such that under the weight of the concrete the final desired form of the shell is achieved. 

The construction process starts with the assembly of the cable net on-site. Then, it is tightened to an initial pre-stressed state. Because of the flexibility of the net and the fabrication and construction tolerances, the built form needs to be corrected to the designed one. A feedback loop is implemented by iteratively measuring the form and applying inputs to the structure. 
After the control of this pre-stressed form of the cable net, the fabric membrane is laid on top and then the concrete is sprayed on. The load distribution of the concrete and reinforcement introduces new uncertainties to the system. Further corrections of the form of the cable net might therefore be required to reach the designed final loaded form.

\subsection{Graph-Theoretical Description}
The cable net is associated with an underlying graph $\mathcal{G}=(\mathcal{N}, \mathcal{E})$. Its $n$ nodes from the node set $\mathcal{N}$ represent the connection points of the cable net, and its $m$ edges from the edge set $\mathcal{E} \subseteq \mathcal{N} \times \mathcal{N}$ correspond to the cable segments of the net. The node set is divided into the two disjoint subsets $\mathcal{N}_F$ of $n_{F}$ free nodes which lie in the interior of the cable net and $\mathcal{N}_B$ of $n_B$ boundary nodes which are attached to the rigid frame. The edge set is composed of two disjoint sets $\mathcal{E}_F$ of $m_{F}$ free edges between free nodes and $\mathcal{E}_B$ of $m_B$ boundary edges which connect the boundary nodes on the rigid frame to free nodes in the interior of the net. It holds that $\mathcal{N} = \mathcal{N}_F \cup \mathcal{N}_B$, $\mathcal{N}_F \cap \mathcal{N}_B = \emptyset$ and $\mathcal{E} = \mathcal{E}_F \cup \mathcal{E}_B$, $\mathcal{E}_F \cap \mathcal{E}_B = \emptyset$. 
By a slight abuse of notation, we use both the edge set and an index set for the edges, i.e., the index $e \in \{1,...,m\}$ or equivalently the index $(s,t) \in \mathcal{E}$ denotes the edge $e$ connecting nodes $s$ and $t$. 

\begin{figure}
\includegraphics[width=0.55\columnwidth]{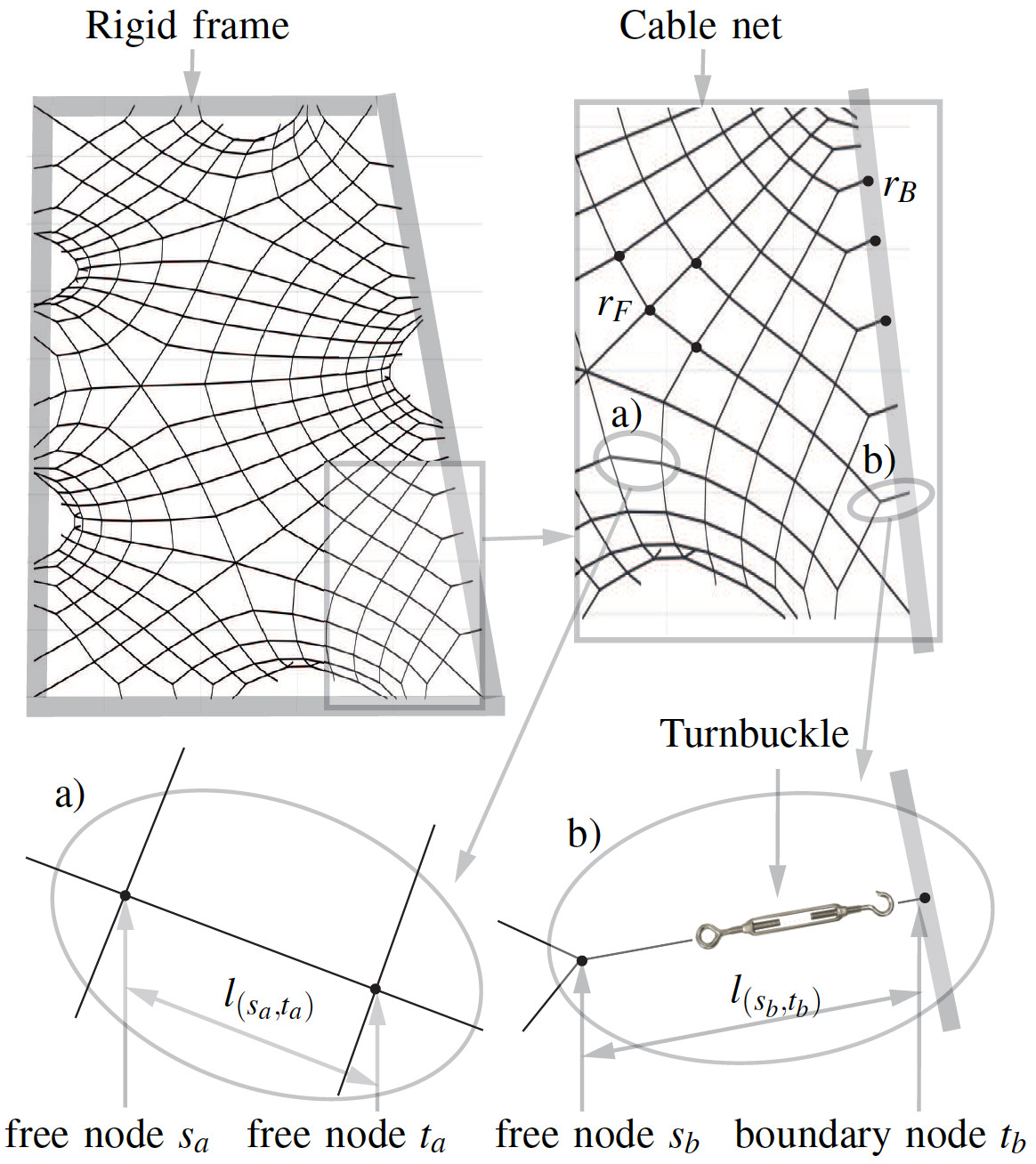}
\caption{Top view of the cable net system, (from \cite{Stuerz2016a})\,. a) Free edge $(s_a,t_a) \in \mathcal{E}_F$ connecting nodes $s_a \in \mathcal{N}_F$ and $t_a \in \mathcal{N}_F$ with edge length $l_{(s_a,t_a)}$.  b) Boundary edge $(s_b,t_b) \in \mathcal{E}_B$ with edge length $l_{(s_b,t_b)}$ connecting the free node $s \in \mathcal{N}_F$ and the boundary node $t \in \mathcal{N}_B$. \label{fig:Skizze}}
\end{figure}

The geometric form of the net is described via the positions of the nodes of the cable net, which we define as 
$${x = [x_F^\trans\,, x_B^\trans]^\trans\,, \quad y = [y_F^\trans\,, y_B^\trans]^\trans\,, \quad z = [z_F^\trans\,, z_B^\trans]^\trans ~ \in \mathbb{R}^{n}\,,}$$ where 
the subscripts distinguish between the free and boundary nodes. For individual nodes, $s \in \{1,...,n\}$, we define the vector of their coordinates as \vspace{-0.0cm} $$r_s = [x_s\, , \, y_s \, , \, z_s]^\trans \in \mathbb{R}^3.$$ \vspace{-0.0cm}
We also define the stacked vectors of the coordinates of all free, all boundary and the collection of all the nodes as 
$$r_F = [r_1^\trans\, , \, ... \, , \, r_{n_F}^\trans]^\trans, \quad r_B = [r_{n_F+1}^\trans\, , \, ... \, , \, r_{n}^\trans]^\trans \text{~and~} r = [r_F^\trans \, , \, r_B^\trans]^\trans.$$ 
We use the term configuration as a synonym for the form of the cable net defined through the nodal position coordinates. 
A top view of the cable net is depicted in Fig.~\ref{fig:Skizze} with the free nodes $s_a$ and $t_a$ connected by the free edge $(s_a,t_a)$ and with the free node $s_b$ connected by the boundary edge $(s_b,t_b)$ to the boundary node $t_b$. 
The actual length of an edge $(s,t)$ is denoted by $l_{(s,t)}$ and is given by the Euclidean distance between its two nodes as 
\begin{align}
l_{(s,t)}:=\|{r}_s-{r}_t\|_2\,.
\label{lst}
\end{align}

\subsection{Parameters and Inputs}
The materials and dimensions of the edges in the cable net are described by the following fixed parameters. The Young's modulus $E$ indicates the relation between stress and strain in the material. 
The constant $EA_{(s,t)}$ is used to denote the product of the Young's modulus $E$ and the cross section area $A$ of the edge $(s,t)$ and defines its elastic properties. 
Other important parameters of the system influencing the forces within the net are the unstressed lengths of the edges, 
which are denoted by 
\begin{equation} \label{eq:l0def}
l_{0} = \begin{bmatrix} l_{0,1}, ... , l_{0,(s,t)}, ... , l_{0,m} \end{bmatrix}^\trans \in \mathbb{R}^m, 
\end{equation}
with $l_{0,(s,t)}$ the parameter of the edge $(s,t)$. In the interior of the net, these parameters are fixed and cannot be changed once the cable net has been constructed. 
The boundary edges connect the cable net to the rigid frame via turnbuckles, as can be seen in Fig.~\ref{fig:Skizze}. These turnbuckles can introduce a defined change in length of the boundary edges, which are used as inputs to the system in order to control the form. 
The possible input vector is thus defined as the vector collecting all the changes in lengths for all boundary edges  
\begin{align}
{u}:=[u_1\, , \, ... \, ,\, u_{(s,t)}\, , \, ... \, , \, u_{m_B}]^\trans \in \mathbb{R}^{m_B}, \quad \forall (s,t) \in \mathcal{E}_B, 
\label{eq:udef}
\end{align} 
with 
$u_{(s,t)}$ being the change in length $l_{0,(s,t)}$ for the boundary edge $(s,t)$. 
Here, the assumption is made that the turnbuckles are not elastic. The unstressed length of edge $(s,t)$ after applying the input $u_{(s,t)}$ is defined as \vspace{-0.2cm}
\begin{equation} \label{eq:l0barst1}
\overline{l}_{0,(s,t)} = l_{0,(s,t)}-u_{(s,t)}\,, ~\forall (s,t) \in \mathcal{E}_B.
\end{equation}
We also use the notation 
\begin{equation} \label{eq:l0barst2} 
\overline{l}_{0,(s,t)} = l_{0,(s,t)}\,, ~\forall (s,t) \in \mathcal{E}_F, 
\end{equation}
if $(s,t)$ is a non-adjustable free edge. In Fig.~\ref{fig:Skizze}, the free edge $(s_a,t_a)$ is of constant unstressed length $\overline{l}_{0,(s_a,t_a)} = {l}_{0,(s_a,t_a)}$, 
and the boundary edge $(s_b,t_b)$ is of adjustable unstressed length $\overline{l}_{0,(s_b,t_b)} = l_{0,(s_b,t_b)}-u_{(s_b,t_b)}$, with possible ${u}_{(s_b,t_b)}\neq 0$. 

With \eqref{lst} and \eqref{eq:l0barst1}, \eqref{eq:l0barst2}, the actual elongation of the edge $(s,t)$ is given by
\begin{equation} \label{eq:Deltal}
\Delta l_{(s,t)} = l_{(s,t)} - \overline{l}_{0,(s,t)}\,. 
\end{equation}
Configurations that have no slack cables avoid sagging of the concrete and are preferred. This means that the edges of the cable net are desired to be in zero or positive tension. 
Therefore, the following constraints can be introduced 
\begin{equation}
\begin{aligned}
g_{(s,t)}({r}_F, r_B,{u}) &:= - \Delta {l_{(s,t)}} \leq 0\,, \quad \forall (s,t) \in \mathcal{E}.
\label{eq:ineqcons}
\end{aligned}
\end{equation}
We summarize these constraints for all edges in the vector 
\begin{equation}
\begin{aligned}
g({r}_F, r_B,{u}) &:= [g_{1}, ... , g_{m}]^\trans \,. 
\label{eq:ineqcons_sum}
\end{aligned}
\end{equation}
However, in the physically built cable net, there may be slack cables because of construction imprecision. These slack edges may or may not be removable by the control, depending on the parameters of the edges. In other words, there may or may not exist a configuration with no slack edges for the given parameters. 

Depending on the construction, there might be physical limitations on the possible change in the boundary edge lengths. Then, input constraints in the form of 
\begin{equation} \label{eq:inputcons}
\begin{aligned}
{u} &\leq {u}_{ub}, \\
-{u} &\leq - {u}_{lb},
\end{aligned}
\end{equation}
might need to be introduced, with ${u}_{ub} \in \mathbb{R}^{m_B}$ and ${u}_{lb} \in \mathbb{R}^{m_B}$ being upper and lower bounds on the possible inputs, 
respectively.

\subsection{Parameter Identification}
The stiffnesses and therefore the forces and the form of the cable net are very sensitive to the unstressed lengths $l_0$, defined in \eqref{eq:l0def}. However, the lengths are subject to fabrication tolerances and uncertainties and are therefore likely to deviate from the nominal values of the design model. 
The cable net might be constructed in its stressed state, such that only the actual stressed lengths of the edges can be measured, but not the unstressed lengths. 

Methods for the parameter identification of the unstressed lengths are proposed in \cite{Stuerz2016b}. They require measurements of different configurations of the cable net, which are obtained by exciting the system with different inputs. The identification is based on the model of the cable net in static equilibria, which is described in the following.

\section{Models of the Static Equilibria}
\label{sec:model}
In order to describe the form of the cable net and its dependency on the inputs, a model of the static equilibria of the cable net is required in terms of the parameters and boundary conditions of the system. Two equivalent formulations were introduced in \cite{Stuerz2016a}. 
One is based on force balances at all free nodes, \cite{Connelly1996}, and leads to implicit nonlinear equations. The other is an energy minimization approach and can be cast as a convex optimization problem. Both models will be used to formulate the control algorithm in Section~\ref{alg:SQP}.

\subsection{Model Assumptions} 
\label{subsec:modelass}
From the first pre-tensioned state of the cable net after its assembly to its final controlled state, a series of equilibria configurations of the cable net are considered. In order to model these equilibria configurations, we define the reduced graph $\bar{\mathcal{G}}$ generated by removing all slack edges. Thus, the graph $\bar{\mathcal{G}}$ is defined as $\bar{\mathcal{G}} = \{ \mathcal{N}, \bar{\mathcal{E}} \}$ with $\bar{\mathcal{E}}$ being the set of all tensioned edges of the cable net. Here, it is assumed that nodes that have only slack adjacent edges have already been removed from the node set $\mathcal{N}$. 
A further assumption is that the cable net is designed in such a way that $\bar{\mathcal{G}}$ does not contain 2-cycles or self-loops. 

We assume that $\bar{\mathcal{G}}$ is known and that it stays constant for the series of equilibria configurations considered during the control, i.e., cables do not change from being slack to being tensioned or vice versa. In the design phase, the desired force distribution in the cable net is designed in such a way that all cables are in positive tension. In practice, checking these conditions in the real system can be done by force measurements, by manual examination, or by measuring whether the actual edge lengths, $l$, are longer in the prestressed state than the initial edge lengths, $l_0$. 

Furthermore, we assume that the parameters of the system, i.e., the material properties $EA$, the unstressed lengths of the edges $l_0$, and the self-weight of the net are known and are constant. These parameters can be identified in experiments. In the following, we do not consider uncertainties on the model parameters. 
The boundary nodes at the rigid frame can be measured and are thus also known. They are considered to stay constant, as they are fixed at the rigid frame.

\subsection{Energy Minimization Approach}
In order to find the static equilibrium of the system, an approach of minimizing its total energy can be taken. In the following, we assume that the elastic tension forces versus elongation function of the edges are linear and increasing.  Note that this can easily be generalized to piecewise linear and increasing functions.  
Under this assumption, for a given fixed vector of inputs $u$ in \eqref{eq:udef}\,, this energy minimization problem is equivalent to a convex 
second-order cone program (SOCP). 
This approach was presented for a similar example in \cite{Lobo1998} and for the cable net application in \cite{Stuerz2016a}. 

The total energy of the cable net, expressed by 
\begin{equation} \label{eq:totalenergy}
\begin{aligned}
V(r,u) =  
-{p}_z {z}_F + \sum_{(s,t) \in \mathcal{E}} \frac{EA_{(s,t)}}{2 {\overline{l}}_{0,(s,t)}} (l_{(s,t)}-\overline{l}_{0,(s,t)})^2, \\
\end{aligned}
\end{equation}
which is the sum of the potential energies of all nodes (first term) and the sum of the elastic energies of all tensioned edges (second term).
The vector ${p}_z$ in \eqref{eq:totalenergy} accounts for point loads due to self-weight and any other loads on all free nodes. 
The problem of finding the equilibrium of the cable net by minimizing this total energy term \eqref{eq:totalenergy} is given as
\begin{align}
& &&r(u) = \underset{r}{\mathrm{argmin}} ~~  V  & && \notag \\
& && ~ \mathrm{s.t.}  ~~~ r_B = \overline{r}_B\,, & && 
\label{eq:minimizeenergy}
\end{align}
where $\overline{r}_B$ are the fixed positions at the rigid frame where the boundary edges are connected. 

For a fixed input vector $u$, i.e. constant $\overline{l}_0$, it is possible to rewrite the problem given in \eqref{eq:minimizeenergy} as a convex optimization problem. We introduce a variable ${v}$ and vector $\w$. The entry ${\w_{(s,t)}}$ of $\w$ for the edge $(s,t)$ is defined as 
\begin{align} \label{eq:wplus}
{\w_{(s,t)}} = \left( \frac{\sqrt{EA_{(s,t)}}}{\sqrt{ {\overline{l}}_{0,(s,t)} }} ~ {(l_{(s,t)}-\overline{l}_{0,(s,t)})}  \right)_+, 
\end{align}
with the notation $\left( \cdot \right)_+=\mathrm{max}(0, \cdot )$. The upper bound ${v}$ on the term $\left\Vert \w \right\Vert_2^2$ fulfills the following hyperbolic constraint 
\begin{align}
\left\Vert \w \right\Vert_2^2 \leq {v} \iff \left\Vert \begin{bmatrix} 2 \w \\ 1- {v} \end{bmatrix} \right\Vert_2 \leq 1+ {v}. \notag
\end{align}
Rewriting the problem in terms of the variables $\w$ and $v$, we find that the coordinates ${r}_F$ of the free nodes in a static equilibrium configuration can be obtained as the minimizers to the following SOCP. 
\begin{equation}\label{SOCP}
   \begin{array}{l l l}
\text{Problem} ~  \mathcal{P}_{\mathrm{minE}}:  & &    \\[0.2cm]  
& \underset{{r}_F, v, \w}{\min} ~  & - {p}_z  {z}_F + \frac{1}{2} {v}    \\[0.3cm]
&  ~\mathrm{s.t.}   & \frac{\sqrt{EA_{(s,t)}}}{\sqrt{ {\overline{l}}_{0,(s,t)} }} ~  \left( \left\Vert  {r}_s  -  {r}_t  \right\Vert_2 - \overline{l}_{0,(s,t)} \right)  \leq \w_{(s,t)} ,   \\[0.3cm]
& & 0 \leq \w_{(s,t)}, \quad \forall ~(s,t) \in \mathcal{E},  \\[0.2cm] 
&  &  {r}_{B} = \overline{r}_{B},    \\[0.2cm]
& & \left\Vert \begin{bmatrix} 2 \w \\ 1- {v} \end{bmatrix} \right\Vert_2 \leq 1+ {v}\,, 
   \end{array}
\end{equation}
for a fixed input vector $u$, i.e. constant lengths $\overline{l}_{0,(s,t)}$. 

Note that the definition of $\w_{(s,t)}$ in \eqref{eq:wplus} allows for only positive tension forces to contribute to the energy $V$ of the system. This is consistent with the model assumption that in the case where $l_{(s,t)} < l_{0,(s,t)}$, the cable is not in compression, but it is a slack cable under zero force. However, the solution of Problem $\mathcal{P}_{\mathrm{minE}}$ does not guarantee that there are no slack cables in the equilibrium state of the cable net. 

\subsection{Force Balance Approach}
The net force at each free node $s$ is the sum over the tension forces of all its adjacent edges $(s,t) \in \bar{\mathcal{E}}_s$, 
where we denote by $\bar{\mathcal{E}}_s$ the set of all adjacent edges of node $s$ which are in tension, i.e., $\forall (s,t), ~$ such that $~ \Delta l_{(s,t)} \geq 0$. 
Thus, the net force at node $s$ is given by 
\begin{equation}
\begin{aligned} \label{eq:hsabs}
 h_{s} 
& = \sum_{(s,t) \in \bar{\mathcal{E}}_s} EA_{(s,t)} \left( \frac{ l_{(s,t)} - \overline{l}_{0,(s,t)}}{{\overline{l}}_{0,(s,t)}} \right) d_{(s,t)}  \\
 & = \sum_{(s,t) \in \bar{\mathcal{E}}_s} EA_{(s,t)} l_{(s,t)} \left( \frac{1}{{\overline{l}}_{0,(s,t)}} -  \frac{1}{l_{(s,t)}}  \right) d_{(s,t)} \\
& = \sum_{(s,t) \in \bar{\mathcal{E}}_s} EA_{(s,t)} \left( r_s - r_t \right) \left( \frac{1}{{\overline{l}}_{0,(s,t)}} - \frac{1}{l_{(s,t)}}  \right),
\end{aligned}
\end{equation}
where $d_{(s,t)} = (r_s - r_t) / l_{(s,t)}$ is the direction vector of the edge $(s,t)$ along which the corresponding force is acting. 

For a fixed input vector $u$, the static equilibrium of the cable net can be characterized by the configuration $r_F$ for which all the net forces at all free nodes are zero, i.e., which is the solution of the equations 
\begin{equation}
\begin{aligned} \label{eq:hsxyz}
 h_{s} 
& = \sum_{(s,t) \in \bar{\mathcal{E}}_s} EA_{(s,t)} \left( \begin{bmatrix} x_s \\ y_s \\ z_s \end{bmatrix} - \begin{bmatrix} x_t \\ y_t \\ z_t \end{bmatrix} \right) \left( \frac{1}{{\overline{l}}_{0,(s,t)}} - \frac{1}{l_{(s,t)}}  \right) = 0, 
\end{aligned}
\end{equation}
for all free nodes $s \in \mathcal{N}_F$ and for all tensioned adjacent edges, $(s,t) \in \bar{\mathcal{E}}_s$ to node $s$. 
Note that summing the forces of only the tensioned adjacent edges in \eqref{eq:hsxyz} 
prevents from accounting for the contribution of slack edges as compression forces. 

The function $h: \mathbb{R}^{3 n} \times \mathbb{R}^{m_B} \mapsto \mathbb{R}^{3 n_F}$ is the vector of all force equilibrium equations for all free nodes, i.e., 
\begin{equation} 
h(r, u) = \begin{bmatrix} h_1^\trans & \hdots & h_{n_F}^\trans \end{bmatrix}^\trans. 
\end{equation} 
For fixed boundary values, $r_B = \bar{r}_B$, we may simplify the notation to 
$h: \mathbb{R}^{3 n_F} \times \mathbb{R}^{m_B} \mapsto \mathbb{R}^{3 n_F}$ with 
\begin{equation} 
\label{eq:h}
h(r_F, u) = \begin{bmatrix} h_1^\trans & \hdots & h_{n_F}^\trans \end{bmatrix}^\trans. 
\end{equation} 

\subsection{Equivalence of the Models} 
We show that for fixed parameters and a fixed input vector, the formulations based on the force equilibria and on the minimum energy each have a unique solution and are thus equivalent. This is important for the convergence guarantees of the control algorithm presented in the next section. 
 
First, we state the following result. 
\begin{myprop} \label{prop:Jacinv}
Under the model assumptions in Section~\ref{subsec:modelass}, and for a constant input $u$, 
the partial Jacobian $\nabla_{r_F} h(r_F,u)$ at an equilibrium configuration $r_F$ is non-singular. 
\end{myprop}
\begin{IEEEproof}
Under the model assumptions in Section~\ref{subsec:modelass}, and for a fixed input $u$, 
it can be shown that there exists a unique equilibrium configuration of the cable net, $r_F$, as there exists a unique minimum of the energy function $V(r,u)$ in \eqref{eq:totalenergy}, \cite{Connelly1982}. 
Therefore, the Hessian of the energy function is positive definite at the equilibrium configuration $r_F$, \cite{Connelly1996}, and the partial Jacobian of the force equilibrium equations, $\nabla_{r_F} h(r_F)$ is equal to the Hessian of the energy function. 
\end{IEEEproof}

We now show that the mapping $u$ to $r_F$ via the force equilibrium equations, i.e., the $r_F$ solving \eqref{eq:hsxyz} for a given $u$, is unique.  
In the following, we make use of the Implicit Function Theorem for the function $h(r_F,u)$ in 
\eqref{eq:hsxyz}, and therefore briefly restate it here. 
\begin{mythe}[Implicit Function Theorem \cite{Nocedal2006}]
\label{the:implicitfunction}
Let $h: \mathbb{R}^{3 n_F} \times \mathbb{R}^{m_B} \mapsto \mathbb{R}^{3 n_F}$ be a function such that 
\begin{enumerate}
\item[(i)] $h(\hat{r}_F,\hat{u})=0$ for some $\hat{r}_F \in \mathbb{R}^{3 n_F}$, 
\item[(ii)] the function $h(\cdot,\cdot)$ is continuously differentiable in some neighborhood of $(\hat{r}_F,\hat{u})$, and
\item[(iii)] $\triangledown_{r_F} h(r_F,u)$ is nonsingular at the point $(r_F,u) = (\hat{r}_F,\hat{u})$. 
\end{enumerate}
Then there exist open sets $\mathcal{N}_{r_F} \subset \mathbb{R}^{3 n_F}$ and $\mathcal{N}_{u} \subset \mathbb{R}^{m_B}$ containing $\hat{r}_F$ and $\hat{u}$, respectively, and a unique continuous function $\frF(u): \mathcal{N}_{u} \mapsto  \mathcal{N}_{r_F}$ such that $\hat{r}_F = \frF(\hat{u})$ and $h(r_F,u) = 0$ for all $u \in \mathcal{N}_{u}$. 
If $h(r_F,u)$ is $p$-times continuously differentiable w.r.t.\ both $r_F$ and $u$ for some $p > 0$, then $\frF(u)$ is also $p$-times continuously differentiable w.r.t.\ $u$, 
and we have 
$$\nabla_u \frF(u) = - [\nabla_{r_F} h(r_F,u) ]^{-1} \nabla_u h(r_F,u), $$ 
for all $u \in \mathcal{N}_u$. 
\end{mythe}
The function $h(r_F,u)$ in \eqref{eq:hsxyz} fulfills (i) because we assume that 
for the given parameters 
and a given $\hat{u}$ and under the 
model assumptions in Section~\ref{subsec:modelass}, 
there exists an equilibrium configuration $\hat{r}_F$. 
The condition (ii) holds due to the function definition of $h(\cdot,\cdot)$ in \eqref{eq:hsxyz}. 
Condition~(iii) is fulfilled because of Proposition~\ref{prop:Jacinv}. 

Due to the Implicit Function Theorem it holds that for known parameters $EA$, $l_0$, fixed boundary points $r_B$, and a fixed input vector $u$, there exists a unique equilibrium state of the system, i.e., there exists an $r_F$ which is the unique solution of the force equilibrium equations \eqref{eq:hsxyz}. 

It remains to show uniqueness of the mapping $u$ to $r_F$, i.e., 
that the minimizer $r_F$ of Problem $\mathcal{P}_{\mathrm{minE}}$ is unique. 
We refer to \cite{Kanno2003} for a detailed proof. 
Therein, the problem of minimizing the total potential energy of a cable net structure is reformulated as an SOCP in standard form, which can easily be dualized. The resulting dual SOCP in standard form can again be reformulated into a problem with a  physical interpretation, which is  the minimization of the total complementary energy of the system. The potential strain energy is defined as the integral of the tension force over the elongation of an edge, whereas the complementary energy is defined as the integral of the elongation of the edge over the tension force. 
In the total complementary energy minimization form it can be shown that the deformation of the cable net from an initial configuration to the equilibrium configuration is unique. It then follows that the primal problem has a unique solution corresponding to the unique minimal total potential energy state. 
With the model assumptions in Section~\ref{subsec:modelass}, 
the cable net system considered here and thus $\mathcal{P}_{\mathrm{minE}}$ has a unique solution $r_F$, which is equal to the unique solution $r_F$ of the force equilibrium equations \eqref{eq:hsxyz}. 

\section{Control Algorithm}
\label{sec:Control}
Both the form of the initial pre-stressed unloaded and concrete loaded cable net 
are likely to deviate from the computed pre-stressed initial form 
and the designed final target form, 
respectively. 
Therefore, a closed-loop construction method is introduced, which means that iteratively, the form of the net is measured and inputs are applied to the system to bring its shape as closely as possible to the target shape. 

After formulating the control problem, we present the iterative control algorithm, which is a variant of SQP, 
and give a short comparison to standard SQP methods. Then, an extension of the algorithm is presented to compute sparse control inputs. 

\subsection{Control Problem Formulation}
An optimal control problem (OCP) is formulated, where the cost function to be minimized is given by the weighted $L_2$-norm of the distance between the measured and the desired coordinates, $r_F$ and $r_F^{\mathrm{des}}$, respectively. 
For the solution of this problem to be a static equilibrium of the cable net, 
constraints are used to represent the static equilibrium conditions. 
The OCP can thus be formulated as 
\begin{equation}\label{eq:optproblem}
   \begin{array}{l l l}
\text{Problem} ~  \Pocp:  & &    \\[0.2cm]  
& \underset{{r}_F, u}{\min} ~  & f_{\mathrm{ocp}}(r_F) = \frac{1}{2} \left\Vert  {r}_F- {r}_{F}^{\mathrm{des}}  \right\Vert_{Q_r}^2  \\[0.2cm]
&  ~\mathrm{s.t.}   & {h}({r}_F, \bar{r}_B, {u}) = {0} \,, 
   \end{array}
\end{equation}
with ${Q}_r$ being a 
weighting matrix. The $3  n_F$ equality constraints in \eqref{eq:hsxyz}, ${h}([{r}_F^\top, \bar{r}_B^\top]^\top, {u}) = {0}$, represent the force balances at all free nodes. 
\begin{myrem}
The $m$ inequality constraints in \eqref{eq:ineqcons}, $g({r}_F, \bar{r}_B, {u}) \leq 0$, which represent non-negative elongations of the edges, can be added. If the problem is feasible, they guarantee the absence of slack cables. 
If there are physical constraints in form of construction limitations on the inputs, then the constraints in \eqref{eq:inputcons} need to be added. 
In the remainder of the paper, we do not consider any constraints on the inputs. 
The underlying assumption is that the design provides all the actuation that is needed for the control task. 
\end{myrem}

\subsection{Control Algorithm: Feasible Variant of SQP}
\label{subsec:SQPGN}
We propose to solve Problem $\Pocp$ 
efficiently by a variant of SQP \cite{Nocedal2006, Boggs1995}, where the iterates are feasible in each iteration, denoted by $\iter$. 
Within each iteration of the algorithm, a Gauss Newton (GN) step generates a descent direction for $\Pocp$, and solving $\mathcal{P}_{\mathrm{minE}}$ ensures that feasible iterates are tested in a line search. 

Problem $\Pocp$ is iteratively approximated as a Quadratic Program (QP) around a sequence of points $p^\iter = [r_F^{\iter \trans}, ~ u^{\iter \trans}]^\trans$ of the current nodal position coordinates and inputs. The cost function of this QP could be obtained by a quadratic approximation of the Lagrangian of $\Pocp$, which would involve the Hessian of the Lagrangian. We instead take the constrained GN approach \cite{Nocedal2006}, \cite{Schwetlick1985}, where we exploit the least-squares structure of the cost function $\focp(r_F)$ in $\Pocp$. 
The GN iteration only uses the first-order term for the approximate Hessian $H$, i.e., 
$$ H = \nabla_{r_F,u} ((r_F^\iter-r_F^\mathrm{des}) Q_r^{\frac{1}{2}})^\trans  \nabla_{r_F,u} ((r_F^\iter-r_F^\mathrm{des}) Q_r^{\frac{1}{2}}) \!\! = \mathrm{diag}(Q_r, 0).$$ 
This approximation has significant computational advantages if the system is large, as no second order information needs to be computed. It is a good approximation if the residuals ${r_F^\iter - r_F^\mathrm{des}}$ are small or nearly affine. 
With ${\Delta p^\iter = [\Delta r_F^{\iter \trans} \!\! , ~ \Delta u^{\iter \trans}]^\trans},$ the QP in iteration $\iter$ is  given by 
\begin{equation}\label{eq:approxQP}
   \begin{array}{l l l}
\text{Problem} ~  \PSQPGN:  & &    \\[0.2cm]  
& \underset{\Delta p^\iter}{\min} ~  & \fSQPGN(\Delta r_F^{\iter}) = \frac{1}{2} \Delta p^{\iter \trans} H \Delta p^{\iter}  + \nabla_{(r_F,u)} \focp  \Delta p^{\iter}   \\[0.2cm]
&  ~\mathrm{s.t.}   & h(r_F^\iter, u^\iter) + \nabla_{(r_F,u)} h(r_F^\iter, u^\iter)^\trans \Delta p^\iter = 0. 
   \end{array}
\end{equation}
The equality 
constraints of $\PSQPGN$ 
are the linearized constraints of $\Pocp$.
With $H$ and $\focp$, $\fSQPGN$ can be simplified to  
\begin{equation} \label{eq:GN_f}
\begin{aligned}
\fSQPGN(\Delta r_F^\iter) 
&= \frac{1}{2} \|  r_F^\iter - r_F^{\mathrm{des}} + \nabla_{(r_F,u)} (r_F^\iter - r_F^{\mathrm{des}})\Delta p^\iter \|_{Q_r}^2 \\
&= \frac{1}{2} \|  r_F^\iter - r_F^{\mathrm{des}} + \Delta r_F^{\iter} \|_{Q_r}^2. 
\end{aligned}
\end{equation} 
We propose the following feasible variant of SQP, for which global convergence is shown in Section~\ref{sec:convergence_GN}. 
Along $\Delta {u}^{\iter}$, which is the partial minimizer of $\PSQPGN$, a line search is performed to find a step length $\alpha^\iter$. 
In the direction of $\Delta {u}^{\iter}$, feasible points of $\Pocp$ (denoted by $r_F(u^{(\iter +1)})$) 
are computed by solving $\mathcal{P}_{\mathrm{minE}}$. 
These feasible points are those $r_F^{(\iter + 1)}$ that together with ${u^{(\iter +1)} = u^\iter + \alpha^\iter \Delta u^\iter}$ fulfill the nonlinear constraints of $\Pocp$.  
The iterates are then given by 
\begin{equation} \label{eq:iterates_SQP_SOCP_LS} 
{p}^{(\iter+1)} = [ r_F(u^\iter + \alpha^\iter \Delta u^\iter)^\trans \!\! , ~\,\, u^{\iter \trans} \!\! + \alpha^\iter \Delta u^{\iter \trans} ]^\trans, 
\end{equation} 
with $\alpha^{\iter}$ being a suitable step length that fulfills the Wolfe conditions in \cite{Nocedal2006} 
\begin{equation}
\begin{aligned}
\label{eq:Wolfe} 
&\focp({r}_F({u}^{\iter + 1})) \,
\leq \, \focp({r}_F({u}^\iter)) + c_1 \nabla_{r_F} \focp({r}_F({u}^\iter))^\trans \Delta\bar{r}_F^\iter, \\[-0.0cm]
&\nabla_{r_F} \focp({r}_F({u}^{\iter + 1}))^\trans \Delta\bar{r}_F^\iter \, \geq \, 
c_2 \nabla_{r_F} \! \focp({r}_F({u}^\iter))^\trans \Delta\bar{r}_F^\iter, 
\end{aligned}
\end{equation}
where $c_1 \in \mathbb{R}$ and $c_2 \in \mathbb{R}$ are constants fulfilling ${0 < c_1 < c_2 < 1}$, and $\Delta\bar{r}_F^\iter$ is given by $$\Delta\bar{r}_F^\iter = r_F(u^\iter + \alpha^\iter \Delta u^\iter) - r_F(u^\iter). $$
The Wolfe conditions guarantee a sufficient decrease and curvature of the cost function $\focp$ at the new iterate. 
Note that only the cost, and not constraint violations, needs to be accounted for because all points are feasible. 
Different line search algorithms have been proposed in the literature. In this work, an inexact line search algorithm with backtracking is chosen, where the step length satisfies the Wolfe conditions. For further details we refer to \cite{Nocedal2006}.  
A suitable step length $\alpha^\iter$ satisfying the Wolfe conditions always exists 
under the mild assumptions that the cost function $\focp : \mathbb{R}^{n_F + m_B} \mapsto \mathbb{R}$ is continuously differentiable and bounded below along the ray 
$\focp(r_F(u^\iter + \alpha^\iter \Delta u^\iter))$ and $\Delta u^\iter$ is a descent direction for $\focp$. 
This result is stated in Lemma~3.1 in \cite{Nocedal2006}. 

The variant of SQP, where feasibility of the nonlinear constraints is maintained in each iteration, results in 
Algorithm~\ref{alg:SQP} below. Global convergence to a KKT point of $\Pocp$ is shown in Section~\ref{sec:convergence_GN}. 
\begin{algorithm}
 \KwData{
 \begin{itemize}
 \item Initial feasible point ${r}_{F}^0, {u}^0$\,;
 \item Target coordinates ${r}_F^{\mathrm{des}}$\,;
 \item Convergence bound  $c_{\mathrm{c}}$\,;
 \end{itemize}}
 \KwResult{
 \begin{itemize}
 \item KKT point $p = [{r}_F({u})^{\trans}, {u}^{\trans}]^\trans$ of $\Pocp$\,, 
 \end{itemize}}
\textbf{Initialization}: 
\begin{itemize}
\item Set ${r}_{F}^{\iter}={r}_{F}^0, \quad {u}^{\iter}={u}^0, \quad \iter=0\,$\,;
\end{itemize}
 \While{$\| {p}^{(\iter+1)} - {p}^{(\iter)} \| \geq c_{c}$}{
 \begin{itemize}
\item Solve $\PSQPGN$ to obain $\Delta {u}^{\iter}$\,;
\item Perform Line Search Algorithm 
to find \\
step length $\alpha^{\iter}$ and next feasible iterate ${p}^{(\iter+1)}$ 
\item Set $\iter \leftarrow \iter+1$\,;
\end{itemize} 
}
\begin{itemize} 
\item Set $p \leftarrow p^{\iter}$\,; 
\end{itemize} 
 \caption{Overall control algorithm to solve $\Pocp$: SQP variant with line search and feasibility of the nonlinear constraints of $\Pocp$ in each iteration.  } 
 \label{alg:SQP}
\end{algorithm}

\subsection{Comparison to Standard SQP}
The standard SQP methods generate iterates 
$p^{(\iter +1)} = p^\iter + \Delta p^\iter,$
with $\Delta p^\iter = [\Delta r_F^{\iter \trans}, ~\Delta u^{\iter \trans}]^\trans$ being the minimizer of $\PSQPGN$. 
Until convergence, all these iterates might be infeasible, 
which has two main drawbacks. 

A line search for guaranteeing global convergence needs to be performed on a merit function, which accounts for both the decrease in the cost as well as the constraint violations, see \cite{Nocedal2006}\,. This requires design parameters that can be difficult to tune. 
In our variant, the cost function can be chosen as the merit function, as all iterates are feasible. 

Generating infeasible iterates itself can be a disadvantage. 
If time is critical the algorithm may need to be terminated before convergence is reached. 
The corresponding control iterates do not correspond to the cost of the current iterate and can also violate force, extension or slackness constraints. 
To overcome this problem, so-called feasibility-perturbed SQP algorithms have been considered in the literature, see for example \cite{Wright2004, Tenny2004}, or \cite{Lawrence2001} in the context of nonlinear MPC. The search direction in each iteration is ``tilted'' to give a next feasible iterate. 
In our SQP variant, we exploit the model information, and efficiently solve the SOCP $\mathcal{P}_{\mathrm{minE}}$ to obtain feasible iterates in each iteration. 

\subsection{An Input Sparsity Approach}
\label{sec:sparsity}
Depending on the construction application and the site conditions, the actuation system might not be fully automated. 
For the experiments on the prototype presented in Section \ref{sec:Experiments} 
actuation was applied manually. 
For large-scale structures with a large number of boundary edges this process can be very time- and labor-intensive. 
Depending on the deviations in the form that need to be corrected, 
it might be efficient to apply inputs to only a (possibly small) subset of the boundary edges rather than to adjust all of them. Simulation results suggest that this might not significantly compromise the performance. 

Motivated by the goal of making the actuation practically feasible, a sparse input vector is computed. 
To do so, an additional term is introduced in the cost function to account for the cardinality of the input vector. 
As proposed in \cite{Candes2007}, we use the 
weighted $l_1$-norm as a convex regularizer for the cardinality. The resulting sparse input vector is therefore denoted by $\ul$ in the following, and the corresponding cable net configuration $r_F(\ul)$ is denoted by $\rFl$. 
The weighted $l_1$-norm is given by $\| W \, \ul \|_{l_1} = \sum_{i} w_i \vert \uli \vert $, with $W$ being a diagonal matrix of the weights $w_i$. 
The cost function is then convex and given by 
\begin{equation} \label{eq:J_l1}
f_{l_1}=f_{\mathrm{ocp}}(\rFl, \ul) +  \gamma \, \| W \, \ul \|_{l_1}, 
\end{equation}
with $\gamma$ a weighting factor. If $\gamma = 0$, the fully actuated solution is achieved, and if $\gamma$ is increased, the solution becomes more and more sparse. 
If the weights $w_i$ 
are chosen to be the inverses of the entries of $\uli$, then this weighted $l_1$-norm is equal to the cardinality of $\ul$. As the entries $\uli$ are not known a priori, the weights cannot be chosen a priori. Therefore, an iterative reweighting scheme is implemented, see 
\cite{Candes2007}, \cite{Asif2013}. In the first iteration the initial problem with $\gamma=0$ is solved. Then, the weights $w_i$ are updated to penalize smaller entries more and more, approximating the cardinality of $\ul$. 

The sparse OCP, denoted by $\mathcal{P}_{\mathrm{ocp},l_1}$, consists of minimizing $f_{l_1}$ subject to the constraints of $\mathcal{P}_{\mathrm{ocp}}$. 
\begin{equation}\label{eq:reweigt_prob}
   \begin{array}{l l l}
\text{Problem} ~  & \mathcal{P}_{\mathrm{ocp},l_1}:   &    \\[0.2cm]  
& \underset{\rFl, \ul}{\min} ~  & f_{l_1} = f_{\mathrm{ocp}}(\rFl, \ul) +  \gamma \, \| W \ul \|_{l_1}   \\[0.2cm]
&  \quad \mathrm{s.t.}   & h(\rFl,\ul) = {0}\,.
   \end{array}
\end{equation}
 In order to solve Problem $\mathcal{P}_{\mathrm{ocp},l_1}$ 
by the novel variant of SQP introduced in Section~\ref{sec:Control} we transform it into the following QP 
\begin{equation}\label{eq:reweigt_prob}
   \begin{array}{l l l}
\text{Problem} ~  & \mathcal{P}_{\mathrm{SQP},l_1}^{\iter}:   &    \\[0.2cm]  
& \underset{\Delta \rFl^{\iter}, \Delta \ul^{\iter}, \beta}{\min} ~  & f_{l_1}^{\mathrm{GN}} = \frac{1}{2} \left\Vert \rFl^{\iter} + \Delta \rFl^{\iter} - r_F^{\mathrm{des}} \right\Vert^2_{{{Q_r}}}  + \gamma ~ (w^{\trans} \beta)   \\[0.2cm]
&  \qquad  \mathrm{s.t.}   & ( \ul^{\iter}+ \Delta \ul^{\iter}) \leq \beta \\[0.2cm]  
&    & - ( \ul^{\iter} + \Delta \ul^{\iter}) \leq \beta \\[0.2cm]  
&    & \nabla_{(r_F,u)} h(\rFl^\iter,\ul^\iter) \Delta p^\iter + h(\rFl^\iter,\ul^\iter) = {0}\,, 
 \end{array}
\end{equation}
with $w$ being the vector of the weights $w_i$. 
Algorithm~\ref{alg:reweight} summarizes the steps for solving $\mathcal{P}_{\mathrm{ocp},l_1}$ by Algorithm~\ref{alg:SQP} together with an iterative reweighting scheme of the $l_1$-norm in the cost function, i.e., with iteratively updating the weights $w$.  
The result is a sparse input vector $\ul$. 
\begin{algorithm}
 \textbf{Data:}
   \begin{itemize}
   \item Parameters $\tau >0$ and $\epsilon>0$\;
   \item Convergence bound  $c_{\mathrm{w}}$\,;
 \end{itemize}
 \textbf{Result:}
   \begin{itemize}
  \item Sparse input $\ul$\;
  \end{itemize}
\textbf{Initialize:}
\begin{itemize}
\item Set $\nu = 0$, $\gamma > 0$, $w^0=0$\;
\item Solve $\mathcal{P}_{\mathrm{ocp},l_1}$ by Algorithm~\ref{alg:SQP} \\
to obtain initial fully actuated solution $u^0$\; 
  \end{itemize}

 \While{$\nu<1$ \textbf{\textup{or}} $\| w^{\nu} - w^{(\nu-1)} \| \geq c_{w}$}{
 \begin{itemize}
 \item Set $\nu \leftarrow \nu+1$\;
 \vspace{0.1cm}
 
 \item Update weights: \quad 
$w_i^{\nu} = \frac{\tau}{|\uli^{(\nu-1)} | + \epsilon}$\;
 \vspace{0.1cm}
 
\item Solve $\mathcal{P}_{\mathrm{ocp},l_1}$ by Algorithm~\ref{alg:SQP} with $w^{\nu}$ \\
to obtain $p^{\nu} =  [\rFl^{\nu \trans}, ~ \ul^{\nu \trans}]^\trans$\;
\end{itemize} }
\begin{itemize}
\item Set $p = p^{\nu}$\;

\end{itemize}
 \label{alg:reweight}
 \caption{\label{alg:reweight} Computation of sparse input vector $\ul$ corresponding to feasible point $p = 
 [\rFl^\trans, ~ \ul^\trans]^\trans,$ which solves the iteratively reweighted problem $\mathcal{P}_{\mathrm{ocp},l_1}$.}
\end{algorithm}

\section{Convergence of the Control Algorithm}
\label{sec:convergence_GN} 
This section is dedicated to the following convergence result. 
\begin{mythe} \label{the:convergenceAlgo}
For the control task from an initial configuration of the cable net to the desired configuration $r_F^{\mathrm{des}}$, let the model assumptions in Section~\ref{subsec:modelass} hold. 
Then, Algorithm~\ref{alg:SQP} converges to a KKT-point of $\Pocp$. 
\end{mythe}

The proof of Theorem~\ref{the:convergenceAlgo} is formally stated at the end of this section. First, we present the following results, which are then used in the proof. 
In Section~\ref{subsec:reform}, 
we reformulate $\Pocp$ by expressing $r_F$ through the (unknown) function $R(u)$, making use of the Implicit Function Theorem. The resulting unconstrained nonlinear problem (denoted by $\Pocpu$) is equal to the merit function in the line search of Algorithm~\ref{alg:SQP}. 
In Section~\ref{subsec:unique}, 
we show that until convergence, $\PSQPGN$ has a unique solution $\Delta u^\iter$ in each iteration. 
In Section~\ref{subsec:GN} 
we show that the unique solution $\Delta u^\iter$ is a descent direction for both $\Pocpu$ and $\Pocp$. 
With these results, proving convergence to a KKT-point of $\Pocp$ is equivalent to proving convergence of the GN-method to a stationary point of $\Pocpu$, which is a standard result from unconstrained optimization. 

\subsection{Reformulation of the Merit Function}
\label{subsec:reform}
First, we reformulate $\Pocp$ into
\begin{equation}
\begin{aligned}
\text{Problem} ~ & && \!\!\!\!\!\!\! \Pocpu: & && \notag \\  
& &&\min_{{u}} ~~~  \focpu(u) = \frac{1}{2} \left\Vert  \frF(u)- {r}_{F}^{\mathrm{des}}  \right\Vert_{Q_r}^2. & &&
\label{eq:optproblem}
\end{aligned}
\end{equation}
In this formulation, the equality constraints ${h}(r_F,{u}) = {0}$  
are implicitly included in the function $\frF(u)$, which exists due to Theorem~\ref{the:implicitfunction}. 
Problems $\Pocp$ and $\Pocpu$ are equivalent in the sense that $(r_F^*, u^*) = \argmin \,\, \Pocp$ if and only if $u^* = \,\, \argmin \Pocpu$. 
This holds because the minimizer of $\Pocp$, $(r_F^*, u^*)$, needs to be feasible, i.e., it needs to fulfill $h(r_F^*,u^*)=0$ and therefore $\frF(u^*) = r_F^*$, and $\Pocp$ and $\Pocpu$ are exact reformulations. 

Therefore, a stationary point of $\Pocpu$ corresponds to a KKT-point of $\Pocp$, as all constraints are feasible. 
In the following, convergence to a stationary point of $\Pocpu$ is shown.

We note that the line search in Algorithm~\ref{alg:SQP} is solving $\mathcal{P}_{\mathrm{minE}}$ in the trial steps and that this is equal to evaluating $R(u)$. Therefore the line search in Algorithm~\ref{alg:SQP} is equal to a line search on the objective function $\focpu(u)$ of $\Pocpu$. 

\subsection{Existence of a Unique Search Direction $\Delta u^\iter$}
\label{subsec:unique}
To prove further results, we show the following properties of the Jacobian of the constraint function. 
\begin{myprop} \label{prop:Jacrank} 
In each iteration, the Jacobian $\nabla_u h(r_F^\iter, u^\iter)$ has singular values uniformly bounded away from zero, i.e., 
\begin{equation} \label{eq:Jacobian_assumption}
\exists ~ \mu > 0 ~\text{such that}~ \| \nabla_u h(r_F^\iter, u^\iter) \, \tilde{u}\|_2 \geq \mu \,\|\tilde{u}\|_2\,, ~\forall \tilde{u} \in \mathbb{R}^{m_B}, 
\end{equation}
for all $(r_F^\iter, u^\iter)$ in a neighborhood of the bounded level set ${\mathcal{L} = \{(r_F^\iter, u^\iter) \vert \focp(r_F^\iter, u^\iter) \leq \focp(r_F^{(0)}, u^{(0)})\}}$, with $(r_F^{(0)}, u^{(0)})$ being the starting point of the iteration.  
\end{myprop}
\begin{IEEEproof}
$\nabla_u h(r_F,u)$ has a structure which can be transformed into three stacked diagonal matrices for a specific ordering of the $n_F$ nodes and the $m_B$ inputs. 
This means that there are $m_B$ linearly independent first-order changes of the force equilibria $h$ under $u$ at the current point. For the cable net system, this holds as each unit input (where exactly one entry of $u$ is non-zero) produces a first-order change in the net force of its nearest interior node in the net in a specific direction, i.e., linearly independent forces in the net are caused by the different inputs.  
\end{IEEEproof}

With Propositions~\ref{prop:Jacinv} and \ref{prop:Jacrank}, we can choose the weighting matrix $Q_r$, such that 
\begin{equation} \label{eq:Qr_cond}
\mathrm{rank} \left( Q_r^{\frac{1}{2}} \left(\nabla_{r_F} h \right)^{-1} \nabla_u h  \right) = m_B.
\end{equation}
Note that any positive definite weighting matrix $Q_r$ trivially fulfills \eqref{eq:Qr_cond}. 

We can now state the following. 

\begin{mylem} \label{lem:SQP}
If $Q_r$ is chosen such that \eqref{eq:Qr_cond} is fulfilled, then in each iteration of Algorithm~\ref{alg:SQP}, $\PSQPGN$ has a unique solution $\Delta p^\iter = [\Delta r_F^{\iter \top}, ~\Delta u^{\iter \top}]^\top$. 
\end{mylem}
\begin{IEEEproof}
To see that Lemma~\ref{lem:SQP} holds, we show that: 
\begin{itemize}
\item[a)]  In each iteration, the equality constraint Jacobian
$\nabla_{(r_F, u)} h(r_F^\iter, u^\iter) \in \mathbb{R}^{3 n_F \times 3 n_F + m_B}$ has full row rank. 
\item[b)] 
The matrix $H$ is positive definite on the tangent space of the constraints, i.e., $\Delta p^{\iter \trans} H \Delta p^\iter > 0,$ $\forall \Delta p^{\iter} \neq 0,$ s.t.  $\nabla_{(r_F,u)} h(r_F,u) \Delta p^\iter = 0$. 
\end{itemize}
The partial Jacobian $\nabla_{r_F} h(r_F) \in \mathbb{R}^{3 n_F \times 3 n_F}$ has full rank at any equilibrium configuration $r_F$, which is given at the feasible iterates, see Proposition~\ref{prop:Jacinv}. Together with Proposition~\ref{prop:Jacrank}, a) holds. 
To see that b) holds, we note that 
$\Delta p^{\iter \trans} H \, \Delta p^\iter = [ \Delta r_F^{\iter \trans} \,\, \Delta u^{\iter \trans} ] \mathrm{diag}(Q_r,0) [ \Delta r_F^{\iter \trans} \,\, \Delta u^{\iter \trans} ]^\trans =   \Delta r_F^{\iter \trans}  Q_r \Delta r_F^{\iter \trans}$.  
For all $\Delta p^{\iter}$, s.t.\ $\nabla_{(r_F,u)} \,\,  h(r_F,u) \,\,  \Delta p^\iter = 0$, this is equal to 
$\Delta u^{\iter \trans} H_u \Delta u^{\iter}$, with $H_u = \nabla_u h^\top \nabla_{r_F} h^{- \top} Q_r \nabla_{r_F} h^{-1} \nabla_u h$. 
As $Q_r$ is chosen s.t. it fulfills \eqref{eq:Qr_cond},  
it holds that $H_u  > 0$ and thus b) holds. 
a) implies the linear independence constraint qualification (LICQ), \cite{Wright2004} and together with 
b), it implies that $\PSQPGN$ has a unique solution. 
\end{IEEEproof}

\subsection{Descent Direction and Convergence Results}
\label{subsec:GN}
So far, we have shown that Problem $\mathcal{P}_{\mathrm{ocp}}$ can be reformulated into the unconstrained Problem $\mathcal{P}_{\mathrm{ocp,u}}$, which is equal to the merit function in the line search of Algorithm~\ref{alg:SQP}. Furthermore, in each iteration of Algorithm~\ref{alg:SQP}, a unique GN-descent-direction $\Delta p^\iter$ for $\PSQPGN$ is found. Therefore, with Section~\ref{subsec:SQPGN}, the GN-descent-direction $\Delta u^\iter$ for $\mathcal{P}_{\mathrm{ocp}}$ is unique. 

In the next step, we will show that the GN descent direction $\Delta u^\iter$ on $\mathcal{P}_{\mathrm{ocp}}$ is equal to the GN direction on $\mathcal{P}_{\mathrm{ocp,u}}$, denoted as $\uPocpu^\iter$ in the following. 
As $\Delta u^\iter_{\mathrm{ocp,u}}$ cannot be computed, $\Delta u^\iter$ is computed in Algorithm~\ref{alg:SQP} instead. However because of $\Delta u^\iter = \uPocpu^\iter$, $\Delta u^\iter$ is also a descent direction for $\Pocpu$. 

Then, the rest of the convergence proof can be reduced to showing convergence of the GN-iteration on $\Pocpu$, which is a standard result from unconstrained optimization. 

First, we note that in each iteration of Algorithm~\ref{alg:SQP}, if $u^\iter$ is not a critical point of $\Pocpu$, then $\Delta u^\iter$ is a descent direction for $\Pocpu$. This is stated in 
\begin{mythe} \label{the:descent}
Unless $\nabla_u \focpu(u^\iter) = 0$, the GN search direction $\Delta u^\iter$ in Algorithm~\ref{alg:SQP} is always a descent direction for $\focpu$, i.e., $\Delta u^{\iter \trans} \nabla_u \focpu(u^\iter) \leq 0$. 
\end{mythe} 
For the proof of Theorem~\ref{the:descent}, the following two arguments, stated in Lemma~\ref{lem:equalGNdirection} and Lemma~\ref{lem:descentu} are needed. 
As $\frF(u)$ is not explicitly known, the GN-direction $\Delta u^\iter$ is computed by solving $\PSQPGN$ in Algorithm~\ref{alg:SQP}, where the linearized force equilibrium equations $h(r_F,u)$ define the relation between $r_F$ and $u$. The first step in the proof of Theorem~\ref{the:descent} is to show that the partial GN search direction $\Delta u^\iter$ of the minimizer of $\PSQPGN$, $\Delta p^\iter = [\Delta r_F^{\iter \trans},~ \Delta u^{\iter \trans}]^\trans$, is equal to the GN direction for $\Pocpu$. 
In order to see this, let us denote the GN direction of $\Pocpu$ by $\Delta \uPocpu^\iter$ in order not to confuse it with the GN direction $\Delta u^\iter$ of $\Pocp$. 
The GN iteration on $\Pocpu$ is defined as
\begin{equation} \label{eq:GN_iter}
{u}^{(\iter + 1)} = u^{\iter} + \alpha^\iter \Delta \uPocpu^\iter, 
\end{equation}
with $\alpha^\iter$ being a step length from a line search satisfying the Wolfe conditions in \eqref{eq:Wolfe}, with $\Delta \uPocpu^\iter$ instead of $\Delta \bar{r}_F^\iter $ and with $\nabla_u$ instead of $\nabla_{r_F}$. 
With the GN approximation of the Hessian 
$\nabla_u^2 \focpu(u^\iter) \approx (Q_r^{\frac{1}{2}}\nabla_{u}\frF(u))^\trans  (Q_r^{\frac{1}{2}} \nabla_{u}\frF(u))$,
the GN search direction $\Delta \uPocpu^\iter$ is obtained by solving 
\begin{equation} \label{eq:GN_searchdirection}
\begin{aligned}
\nabla_{u}\frF(u^\iter)^\trans Q_r \nabla_{u}\frF(u^\iter) \Delta \uPocpu^\iter 
&= - \nabla_u \focpu(u^\iter) \\
&= - \nabla_{u}\frF(u^\iter)^\trans  Q_r (\frF(u^\iter)-r_F^{\mathrm{des}}). 
\end{aligned}
\end{equation}
We can now state the following. 
\begin{mylem} \label{lem:equalGNdirection}
In each iteration $\iter$, the GN search direction $\Delta \uPocpu^\iter$ for $\Pocpu$ is equal to the partial GN search direction $\Delta u^\iter$ of $\PSQPGN$ in Algorithm~\ref{alg:SQP}. 
\end{mylem}
The proof of Lemma~\ref{lem:equalGNdirection} is given in the Appendix. 

The second step in the proof of Theorem~\ref{the:descent} is to show that the GN-direction $\Delta \uPocpu$ is a descent direction for $\Pocpu$. This is a standard result, stated in the following. 
\begin{mylem}[\cite{Nocedal2006}] \label{lem:descentu}
Unless $\nabla_u \focpu(u^\iter) = 0$, the GN search direction $\Delta \uPocpu^\iter$ is always a descent direction for $\focpu$, i.e., $\Delta \uPocpu^{\iter \trans} \nabla_u \focpu(u^\iter) \leq 0$. 
\end{mylem} 

The equality of $\Delta u^\iter = \Delta \uPocpu^\iter$ in Lemma~\ref{lem:equalGNdirection} together with the result in Lemma~\ref{lem:descentu} complete the proof of Theorem~\ref{the:descent}, that $\Delta u^\iter$ in Algorithm~\ref{alg:SQP} is always a descent direction of $\focpu$. 

For the proof of Theorem~\ref{the:convergenceAlgo}, it remains to be shown that in each iteration, sufficient decrease in the cost is achieved. This follows as a standard result for the line search on the GN method, see Theorem~10.1 in \cite{Nocedal2006}. 

Summarizing the previous results, we can now complete the proof of Theorem~\ref{the:convergenceAlgo} for convergence of Algorithm~\ref{alg:SQP} to a KKT-point of $\Pocp$. 

\begin{IEEEproof}[Proof of Theorem~\ref{the:convergenceAlgo}]
The results of Lemma~\ref{lem:SQP}, Theorem~\ref{the:descent} and 
Theorem~10.1 in \cite{Nocedal2006} 
prove the convergence of Algorithm~\ref{alg:SQP} to a critical point $u$ of $\Pocpu$, which corresponds to $p = [r_F(u)^\trans,~ u^\trans]^\trans$ being a KKT-point of $\Pocp$, because of constraint satisfaction. Thus, the proof of Theorem~\ref{the:convergenceAlgo} is complete.
\end{IEEEproof}

\section{Experimental Results}
\label{sec:Experiments}
This section presents experimental results conducted on a cable net prototype. 
The main goals of the experiments are the validation of the cable net model, as well as the evaluation of the control performance. 

\subsection{Prototype}
\begin{figure}
\begin{center}
\includegraphics[width=0.4\columnwidth]{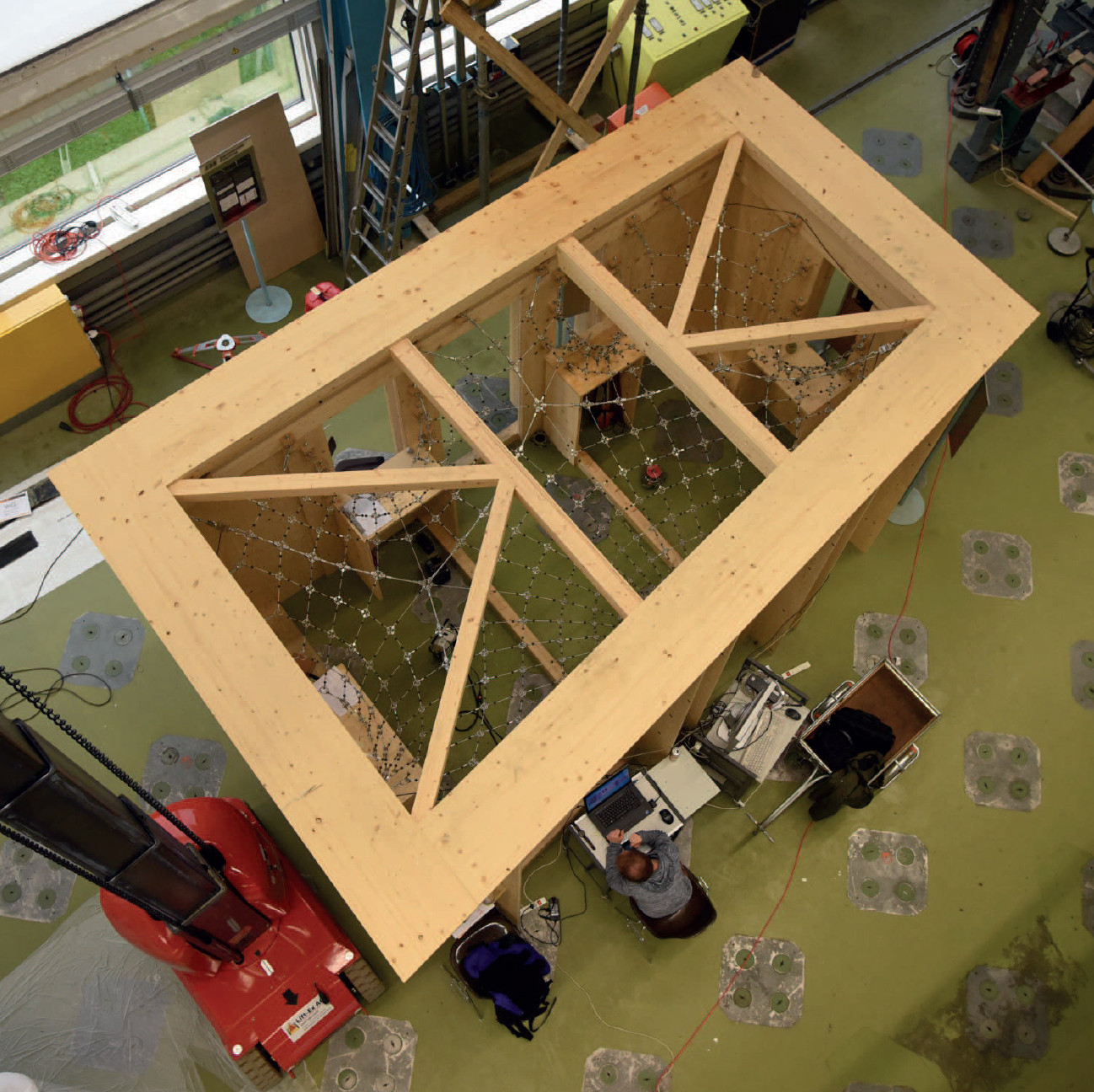}
\end{center}
\caption{\label{fig:exp:topview} Topview of the prototype with stiff wooden frame and cable net structure.}
\end{figure}

\begin{figure}
\begin{center}
\includegraphics[width=0.5\columnwidth]{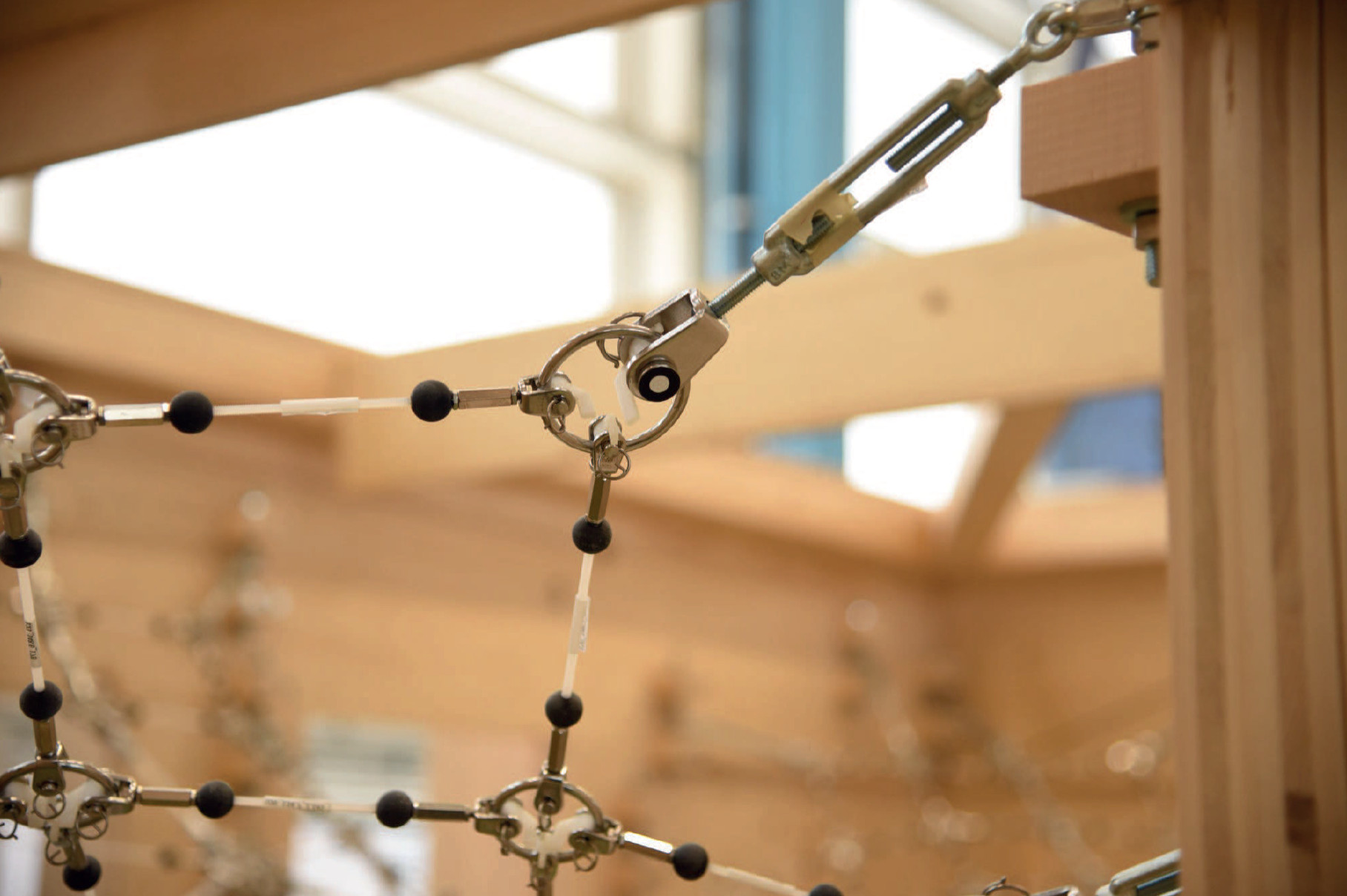}
\end{center}
\caption{\label{fig:exp:turnbuckle} Turnbuckle connecting the cable net structure to the boundary node at the wooden frame and realization of a node in the net as ring construction with attached rods and with markers for the image-based measurements.}
\end{figure}

The experimental prototype is based on the design of the HiLo Roof, which will be built on the NEST building, as described in Section~\ref{sec:intro}. The model is on a scale of 1:4 and therefore its dimensions are approximately $\SI{4.5}{\m} \times \SI{2.5}{\m} \times \SI{2}{\m}$. The rigid frame supporting the pre-stressed net structure is built as a timber housing. A top view of this prototype is shown in Figure~\ref{fig:exp:topview}. The net is realized using plastic and metal rods rather than wire cables. This however is compatible with our model assumptions as the configurations are chosen to be states where all edges are in tension. In simulation and from measurements, this can be verified by computing $\Delta l$. 
The edges are connected via steel connectors to steel rings representing the nodes of the net. The net has a total of $n=295$ nodes, and $m=606$ edges, including $m_B=75$ boundary edges, which are connected via turnbuckles to the anchored boundary points on the wooden frame. 
The realization of the net and the connection to the frame can be seen in Figure~\ref{fig:exp:turnbuckle}. A more detailed description of the prototype system can be found in \cite{Liew2017}.

\begin{figure}
\begin{center}
\includegraphics[width=0.5\columnwidth]{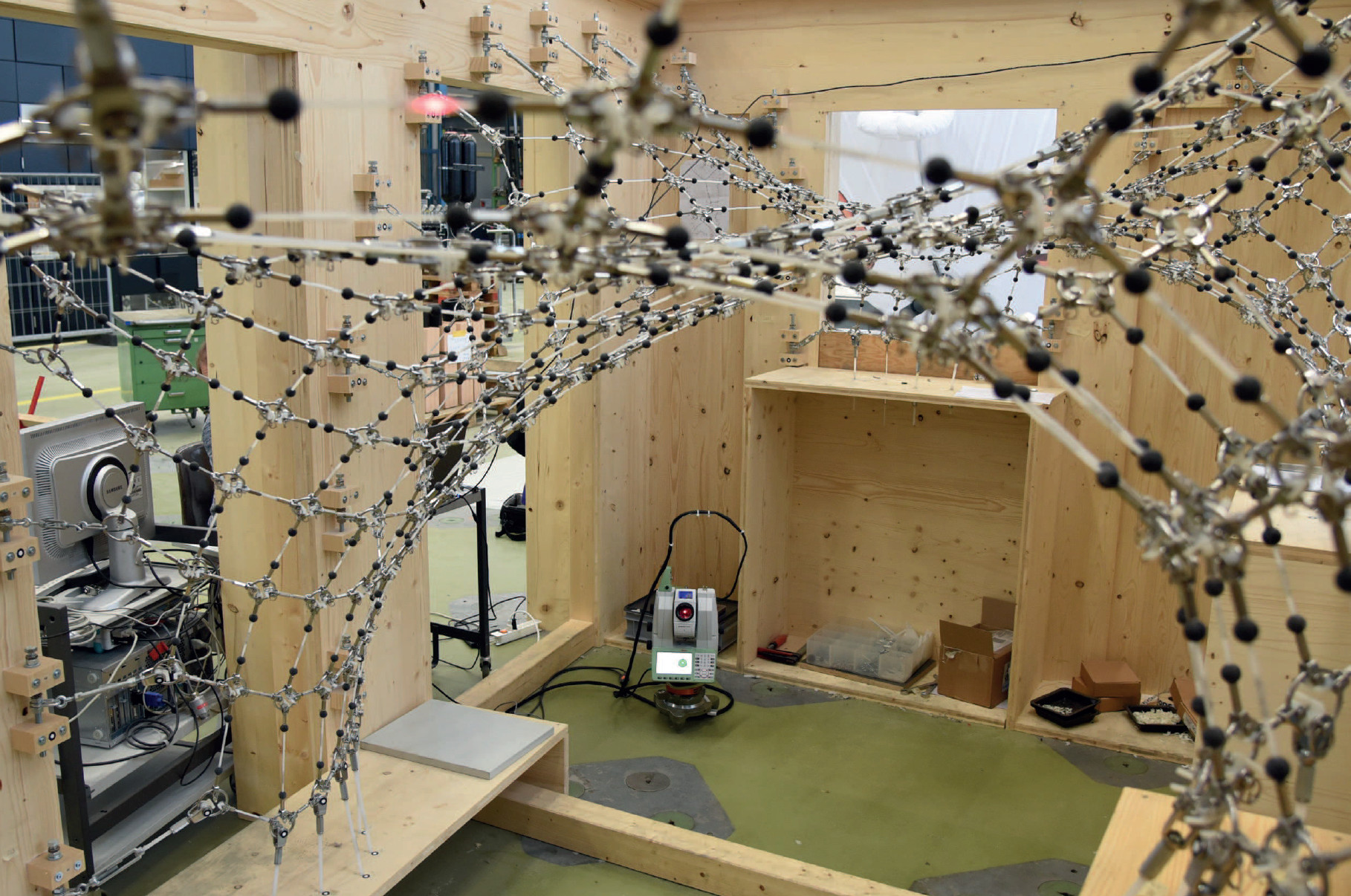}
\end{center}
\caption{\label{fig:exp:theodolite} Theodolite station underneath the net structure pointing towards a marker.}
\end{figure}

\subsection{Measurements}
During the experiments, different cable net configurations are measured, and the $k$-th configuration is denoted by a superscript $(k)$, i.e., $r_F^{(k)}$, or by $r^{(k)}$. 

Whereas the ideal model of the net consists of lines and ideal intersection points, the nodes of the prototype are realized as ring elements (shown in Figure~\ref{fig:exp:turnbuckle}). This makes the estimation of the nodal positions of the net more challenging. The nodes of the ideal model in Section~\ref{sec:model} are defined as the points, where the forces of the adjacent edges balance. 
For the prototype, these points lie close to the center of the ring elements. 

For the experiments, the measurement method is based on image processing. 
The prototype is equipped with black spherical markers, which can be seen in Figure~\ref{fig:exp:turnbuckle}. 
Because of the construction, they cannot be directly attached to the nodes in the center of the rings. Therefore, instead of directly measuring the nodal positions, the measured marker positions will be used to estimate the nodal positions. 

In order to measure the marker positions, the spatial directions from known reference positions to the markers are measured. 
Based on these measurements, the positions of the markers are calculated via triangulation. 
For the measurements of the directions to the markers the vision-based theodolite system QDaedalus, see \cite{Guillaume2008, Guillaume2012, Guillaume2016, Burki2010}, is used in a semi-automated way. In Figure~\ref{fig:exp:theodolite}, the theodolite station which is used to measure the marker positions is shown under the prototype cable net structure.  
For more details about the measurement method, we refer to \cite{Liew2017}. 

The measured marker positions have submillimeter accuracy. We estimate the nodal positions by interpolating the marker positions around each node. The accuracy of this interpolation depends on the locations and number of markers around each node. Near the boundary of the net, the estimated nodal positions were found to be less accurate due to a small number of markers and non-uniform marker placement around the nodes. 

\subsection{Experiments}
The experiments on the prototype are conducted in two phases. In the first phase, the model parameters are determined and the model is validated. In the second phase, control tasks are performed, and the control performance is evaluated. A detailed description of all measured configurations of the prototype net system is presented in \cite{Liew2017}.

\subsubsection{Parameter Identification and Model Validation}
Methods for the identification of the unstressed lengths of the edges are proposed in \cite{Stuerz2016b}. 
The measurement procedure gives precise marker positions. Based on interpolation, we can estimate the nodal displacements precisely, but the estimated absolute nodal positions are not very accurate. 
Therefore, we use a simplified approach for determining the $l_0$ values. A simple model is chosen, where the same material properties are assumed for all edges. 
The first configuration of the net is chosen such that its stress state is approximately uniform. The forces of the edges are at the lower range of possible forces for linear-elastic behavior of the material, however still sufficient to avoid slack edges in the net. 
This was ensured by force measurements at the upper corners of the net and by manual inspection of the edges. 
From the measured nodal coordinates of this initial configuration, denoted by $r^{(0)}$, the parameters of the  unstressed lengths of the edges are determined to be 
\begin{equation}
\begin{aligned} \label{eq:par}
l_{0,(s,t)} &= 0.990 ~ l_{(s,t)}^{(0)}, \quad \text{if~} (s,t) \text{~ is an elastic edge (plastic rod),}  \\
l_{0,(s,t)} &= 0.999 ~ l_{(s,t)}^{(0)}, \quad \text{if~} (s,t) \text{~ is a stiff edge (metal rod),} 
\end{aligned}
\end{equation}
where $l_{(s,t)}^{(0)} = \| r_{s}^{(0)} - r_{t}^{(0)} \|_2$ is the actual measured length of the edge $(s,t)$ in the initial configuration $0$. The stiff edges of the prototype are made of metal because of construction constraints. 
The choice of the scalings for $l_0$ in \eqref{eq:par} is based on the following relation, which holds for the plastic rods, 
\begin{align*}
f_{(s,t)}^{\mathrm{elast}}&=EA_{(s,t)} 
\frac{\Delta l_{(s,t)}}{l_{0,(s,t)}}  \iff \\
 \frac{\Delta l_{(s,t)}}{l_{0,(s,t)}} &= \frac{f_{(s,t)}^{\mathrm{elast}}}{EA_{(s,t)} 
 }  = \frac{\SI{150}{\N}}{1.65 \times 10^9 \frac{\SI{}{\N}}{\SI{}{\m^2}} \pi (\SI{32.4}{\mm})^2}    
 \approx 0.01.   
\end{align*}
The strain $\Delta l / l_0$ of the plastic rod edges is therefore $\approx 1 \%$. 
Based on the assumption of a uniform stress state in the measured configuration, also a uniform strain is assumed for all plastic edges. 
We now make a further simplifying assumption. 
As the plastic rods constitute the majority of the cable edges, they dominate the model behavior. 

For the metal edges, the Young's modulus is higher by approximately a factor of $100$. 
Precisely accounting for this would lead to a model for which the derivatives of the edge forces would have values in significantly different ranges and the OCP would become numerically harder to solve. 
Therefore, the parameters of $l_{0,(s,t)}$ for the stiff edges are chosen as in \eqref{eq:par}. For comparison, the model with higher values of stiffness for the metal rods was used in a simulation study with a small step size and a large number of iterations to convergence. It was confirmed that the simplified model assumptions do not compromise the precision of the simulation results. 

After this initial measurement for determining the unstressed lengths, the subsequent experiments in the first phase are used to validate the model behavior. 
To this end, the nodal positions are estimated from the measured marker positions of the configurations.  
As references, we compute (simulate) the nodal positions based on the model, the estimated parameters from \eqref{eq:par}, the measured boundary coordinates and the inputs of the configurations. 
Then, we compare the measurement-based estimated configurations to the model-based simulated configurations. 

This comparison shows a good match between the simulated model behavior and the behavior of the experimental prototype in terms 
of the displacement of the nodes under the applied inputs. 
However, the match between the simulated and estimated absolute position coordinates is not very accurate for some nodes and has a large variation over the net. A reason for this are the inaccuracies in the nodal position estimates introduced by interpolating the measured marker positions from a small number of markers, and from non-uniformly placed markers around the nodes.  
This appears especially at the boundary regions of the net due to construction limitations. 

To reduce the effect of this estimation error in the absolute nodal positions, the weighting matrix $Q_r$ in the control algorithm, Algorithm \ref{alg:SQP}, is chosen such that it gives more weight to the coordinates that more precisely match the simulated nodal coordinates. This leads to control inputs that correct for the control error rather than correcting for the estimation errors. 

\subsubsection{Control of the Nodal Positions}
In the second phase of the experiments, the goal is to evaluate the control performance. 
We use the superscripts $(\mathrm{ini})$, $(\mathrm{des})$, and $(\mathrm{con})$ to denote the initial, desired and controlled configurations, i.e., $r^{(\mathrm{ini})}$, $r^{(\mathrm{des})}$, and $r^{(\mathrm{con})}$, respectively. 
As before, the controlled configuration resulting from a sparse input vector $\ul$ is denoted by $\rl^{(\mathrm{con})}$. 
The control task is to achieve a given desired target configuration, $r^{(\mathrm{des})}$, starting from an initial configuration, $r^{(\mathrm{ini})}$. 
We present the results of one of several control experiments, as it is representative for the observed control performance. 

\begin{figure}
\begin{center}
\includegraphics[width=0.5\columnwidth]{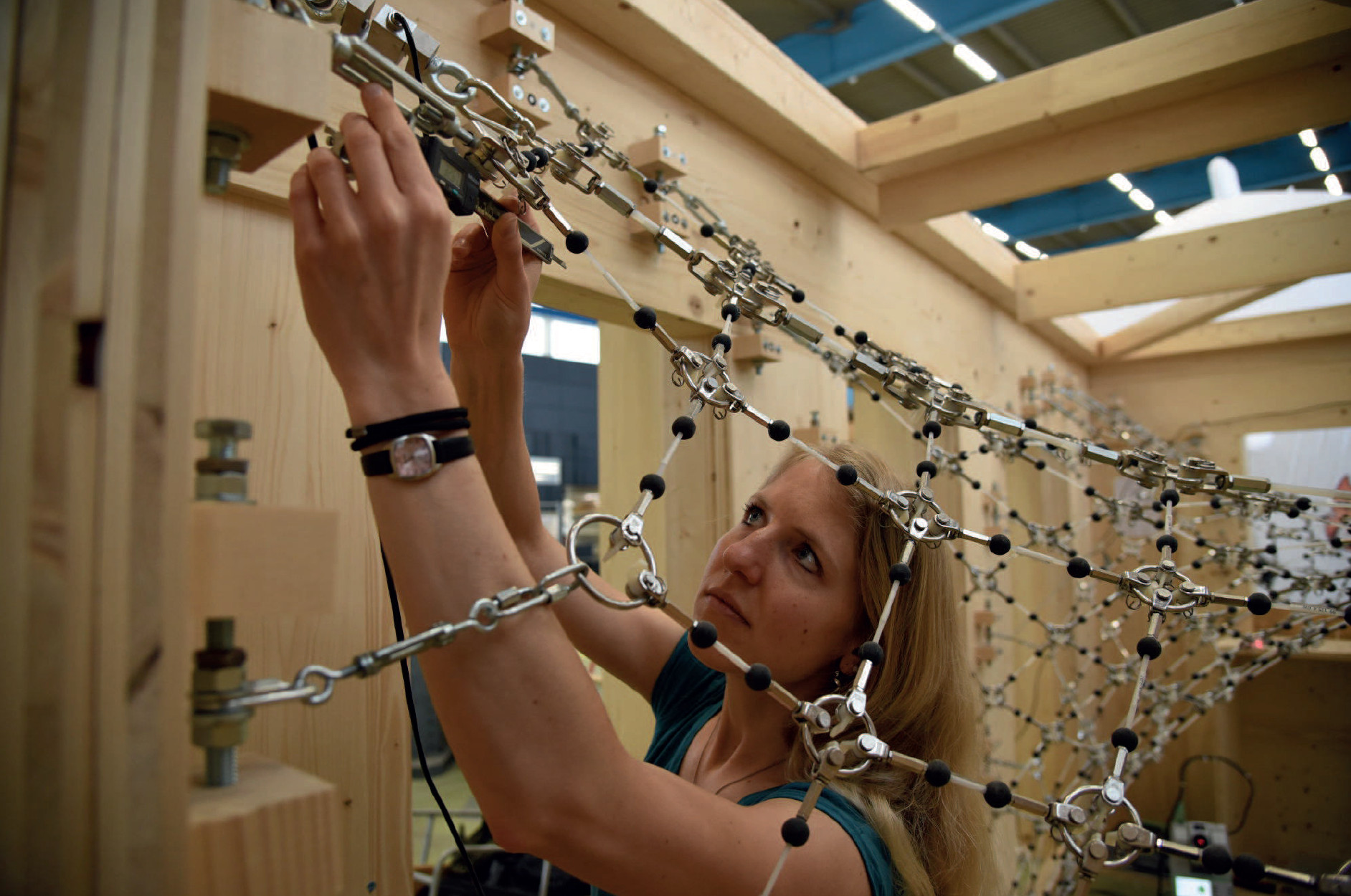}
\end{center}
\caption{\label{fig:exp:control_actions} Manual application of control actions: Measuring the changes in lengths of the boundary edges by callipers and adjusting the turnbuckles to apply the computed control inputs.}
\end{figure}

Figure~\ref{fig:exp:control_actions} shows the manual process of applying the computed control inputs to the prototype system and a motivation for computing sparse input vectors. The turnbuckles are manually actuated to adjust the lengths of the corresponding boundary edges, while measuring the change in lengths via callipers. 

\begin{figure}
\includegraphics[width=0.6\columnwidth]{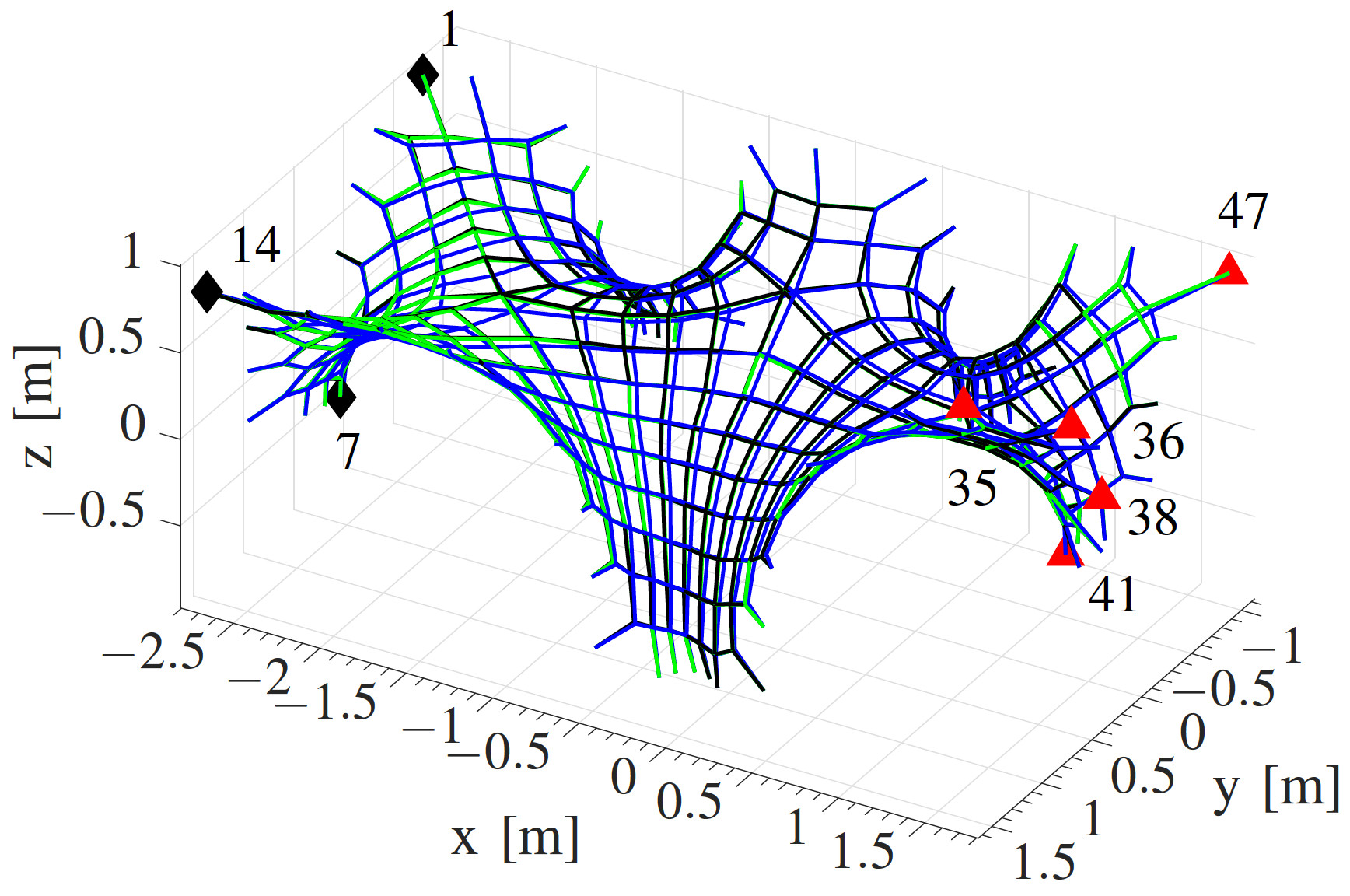}
\caption{\label{fig:Proto_targ_pert_cont} Measured configurations: \textcolor{blue}{\textbf{---} (blue)} Initial starting configuration $r^{\mathrm{(ini)}}$, \textbf{---} (black) Desired configuration $r^{\mathrm{(des)}}$, 
 \textcolor{ForestGreen}{\textbf{---} (green)} Achieved controlled configuration $\rl^{(\mathrm{con})}$. 
 Displacements of configurations $r^{(\mathrm{des})}$ and $\rl^{(\mathrm{con})}$ relative to starting configuration $r^{(\mathrm{ini})}$ are shown scaled by a factor of five for better visualization. 
\textcolor{red}{\textbf{$\blacktriangle$}} Lengthened edges, 
\textcolor{black}{\textbf{$\blacklozenge$}} Shortened edges. The displayed numbers of the actuated edges correspond to a consecutive edge numbering. }
\end{figure}

In Figure~\ref{fig:Proto_targ_pert_cont}, three measured configurations of the 1:4 HiLo roof prototype are shown, which correspond to one control experiment. Blue depicts the initial configuration denoted by $r^{(\mathrm{ini})}$. Black shows the desired target configuration, $r^{(\mathrm{des})}$, and green shows the resulting controlled configuration, $\rl^{(\mathrm{con})}$, after applying the computed sparse inputs. The red triangles \textcolor{red}{\textbf{$\blacktriangle$}}, and black diamonds \textcolor{black}{\textbf{$\blacklozenge$}},  show the actuated boundary edges, which are lengthened and shortened, respectively. 

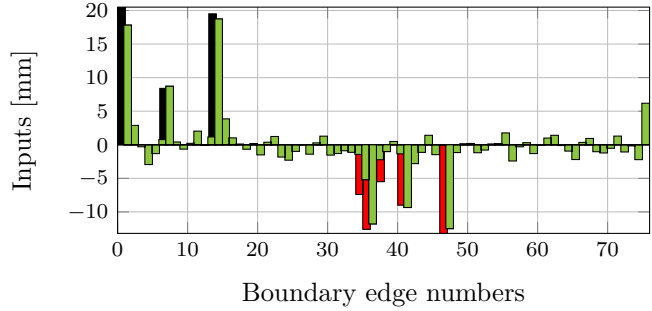
\begin{figure}[h!]
\centering
%
%
\begin{tikzpicture}

\begin{axis}[%
width=7cm,
height=3cm,
at={(1.535433in,0.519685in)},
scale only axis,
area legend,
xmin=0,
xmax=76,
xlabel={Boundary edge numbers},
xmajorgrids,
xtick={0,10,20,30,40,50,60,70},
ytick={-10,-5,0,5,10,15,20},
ymin=-13.2,
ymax=20.5,
ylabel={Inputs [mm]},
ymajorgrids,
ticklabel style = {font=\footnotesize},
]
\addplot[ybar,bar width=0.1cm,bar shift=-0.04cm,draw=black,fill=black] plot table[row sep=crcr] {%
1	20.5\\
2	0\\
3	0\\
4	0\\
5	0\\
6	0\\
7	8.4\\
8	0\\
9	0\\
10	0\\
11	0\\
12	0\\
13	0\\
14	19.5\\
15	0\\
16	0\\
17	0\\
18	0\\
19	0\\
20	0\\
21	0\\
22	0\\
23	0\\
24	0\\
25	0\\
26	0\\
27	0\\
28	0\\
29	0\\
30	0\\
31	0\\
32	0\\
33	0\\
34	0\\
35	0\\
36	0\\
37	0\\
38	0\\
39	0\\
40	0\\
41	0\\
42	0\\
43	0\\
44	0\\
45	0\\
46	0\\
47	0\\
48	0\\
49	0\\
50	0\\
51	0\\
52	0\\
53	0\\
54	0\\
55	0\\
56	0\\
57	0\\
58	0\\
59	0\\
60	0\\
61	0\\
62	0\\
63	0\\
64	0\\
65	0\\
66	0\\
67	0\\
68	0\\
69	0\\
70	0\\
71	0\\
72	0\\
73	0\\
74	0\\
75	0\\
}; \label{refpos}
\addplot[ybar,bar width=0.1cm,bar shift=-0.04cm,draw=black,fill=red] plot table[row sep=crcr] {%
1	0\\
2	0\\
3	0\\
4	0\\
5	0\\
6	0\\
7	0\\
8	0\\
9	0\\
10	0\\
11	0\\
12	0\\
13	0\\
14	0\\
15	0\\
16	0\\
17	0\\
18	0\\
19	0\\
20	0\\
21	0\\
22	0\\
23	0\\
24	0\\
25	0\\
26	0\\
27	0\\
28	0\\
29	0\\
30	0\\
31	0\\
32	0\\
33	0\\
34	0\\
35	-7.4\\
36	-12.6\\
37	0\\
38	-5.5\\
39	0\\
40	0\\
41	-9\\
42	0\\
43	0\\
44	0\\
45	0\\
46	0\\
47	-13.2\\
48	0\\
49	0\\
50	0\\
51	0\\
52	0\\
53	0\\
54	0\\
55	0\\
56	0\\
57	0\\
58	0\\
59	0\\
60	0\\
61	0\\
62	0\\
63	0\\
64	0\\
65	0\\
66	0\\
67	0\\
68	0\\
69	0\\
70	0\\
71	0\\
72	0\\
73	0\\
74	0\\
75	0\\
}; \label{refneg}
\addplot[ybar,bar width=0.1cm,bar shift=0.04cm,draw=black,fill=LimeGreen] plot table[row sep=crcr] {%
1	17.8221265868795\\
2	2.88633471615637\\
3	-0.303690803211095\\
4	-2.94797057880084\\
5	-1.31668664377608\\
6	0.76938220713791\\
7	8.71387166601433\\
8	0.4127995635026\\
9	-0.651908676274597\\
10	0.212016769457264\\
11	2.03338360745727\\
12	-0.0529962518370146\\
13	1.16733549027799\\
14	18.7383835771347\\
15	3.85787761131199\\
16	1.0275260084028\\
17	0.111480818053833\\
18	-0.668308882397736\\
19	0.177406254614259\\
20	-1.50969386194349\\
21	0.380799157447818\\
22	1.23463873480087\\
23	-1.83073465849404\\
24	-2.29837468598877\\
25	-0.982107820284356\\
26	0.00705384640180048\\
27	-1.39987296616536\\
28	0.282548309263559\\
29	1.2822676885956\\
30	-1.54114583006158\\
31	-1.30225401575577\\
32	-0.891120119665366\\
33	-1.12666275517748\\
34	-1.43826423788392\\
35	-5.20192143469544\\
36	-11.7899708011329\\
37	-2.23231898614124\\
38	-1.01148549861266\\
39	0.493929997157961\\
40	-1.34873868163131\\
41	-9.34574245413208\\
42	-2.81171425109248\\
43	-1.13220087974524\\
44	1.41195195779411\\
45	-1.47319037896695\\
46	-0.0226240074902078\\
47	-12.4941155885312\\
48	-1.16167900138323\\
49	0.154091789564874\\
50	0.192971820632532\\
51	-1.18855680427731\\
52	-0.790250459559605\\
53	0.125680364338358\\
54	0.17412189843047\\
55	1.76434557165082\\
56	-2.41999375856993\\
57	-0.30276786544737\\
58	0.312625571639424\\
59	-1.31517120919614\\
60	-0.0389163984896655\\
61	0.998709081097781\\
62	1.40909216585231\\
63	-0.113091245709059\\
64	-0.944595053717183\\
65	-2.20177438338261\\
66	0.354324086237068\\
67	0.940553406942901\\
68	-1.03417204359281\\
69	-1.21920186558534\\
70	-0.522368319340052\\
71	1.29275243054345\\
72	-1.07394851844922\\
73	-0.168260263031386\\
74	-2.21833165514724\\
75	6.1880140184943\\
}; \label{dense}
\addplot[ybar,bar width=0.1cm,bar shift=0.04cm,draw=black,fill=LimeGreen] plot table[row sep=crcr] {%
1	17.8221265868795\\
2	0\\
3	0\\
4	0\\
5	0\\
6	0\\
7	8.71387166601433\\
8	0\\
9	0\\
10	0\\
11	0\\
12	0\\
13	0\\
14	18.7383835771347\\
15	0\\
16	0\\
17	0\\
18	0\\
19	0\\
20	0\\
21	0\\
22	0\\
23	0\\
24	0\\
25	0\\
26	0\\
27	0\\
28	0\\
29	0\\
30	0\\
31	0\\
32	0\\
33	0\\
34	0\\
35	0\\
36	0\\
37	0\\
38	0\\
39	0\\
40	0\\
41	0\\
42	0\\
43	0\\
44	0\\
45	0\\
46	0\\
47	0\\
48	0\\
49	0\\
50	0\\
51	0\\
52	0\\
53	0\\
54	0\\
55	0\\
56	0\\
57	0\\
58	0\\
59	0\\
60	0\\
61	0\\
62	0\\
63	0\\
64	0\\
65	0\\
66	0\\
67	0\\
68	0\\
69	0\\
70	0\\
71	0\\
72	0\\
73	0\\
74	0\\
75	0\\
};
\addplot[ybar,bar width=0.1cm,bar shift=0.04cm,draw=black,fill=LimeGreen] plot table[row sep=crcr] {%
1	0\\
2	0\\
3	0\\
4	0\\
5	0\\
6	0\\
7	0\\
8	0\\
9	0\\
10	0\\
11	0\\
12	0\\
13	0\\
14	0\\
15	0\\
16	0\\
17	0\\
18	0\\
19	0\\
20	0\\
21	0\\
22	0\\
23	0\\
24	0\\
25	0\\
26	0\\
27	0\\
28	0\\
29	0\\
30	0\\
31	0\\
32	0\\
33	0\\
34	0\\
35	-5.20192143469544\\
36	-11.7899708011329\\
37	0\\
38	-1.01148549861266\\
39	0\\
40	0\\
41	-9.34574245413208\\
42	0\\
43	0\\
44	0\\
45	0\\
46	0\\
47	-12.4941155885312\\
48	0\\
49	0\\
50	0\\
51	0\\
52	0\\
53	0\\
54	0\\
55	0\\
56	0\\
57	0\\
58	0\\
59	0\\
60	0\\
61	0\\
62	0\\
63	0\\
64	0\\
65	0\\
66	0\\
67	0\\
68	0\\
69	0\\
70	0\\
71	0\\
72	0\\
73	0\\
74	0\\
75	0\\
};
\end{axis}
\end{tikzpicture}%

\caption{\label{fig:bar_ref_dense} 
\ref{refpos} Reference inputs  for shortening and 
\ref{refneg} for lengthening to steer the initial perturbed configuration  $r^{\mathrm{(ini)}}$ to the desired target configuration $r^{\mathrm{(des)}}$,  
\ref{dense} Computed fully actuated control inputs from Algorithm \ref{alg:SQP}.}
\end{figure}

In order to be able to evaluate the control performance, 
the target configuration $r^{(\mathrm{des})}$ is defined by measuring an actual configuration. This has the advantage of knowing that the target is achievable, and also specifies the inputs required to achieve it. 
Furthermore, we know that the stress states corresponding to this configuration lies within the range of the linear material behavior of the edges, and that no slack edges are present. The control experiment is then started from an initial configuration, $r^{\mathrm{(ini)}}$, which is achieved by perturbing the target configuration, $r^{\mathrm{(des)}}$. The inputs that would then result in the target configuration are depicted in Figure~\ref{fig:bar_ref_dense} as a reference. Their input locations correspond to the ones shown in Figure~\ref{fig:Proto_targ_pert_cont} by red triangles \textcolor{red}{\textbf{$\blacktriangle$}}, and black diamonds \textcolor{black}{\textbf{$\blacklozenge$}}. 

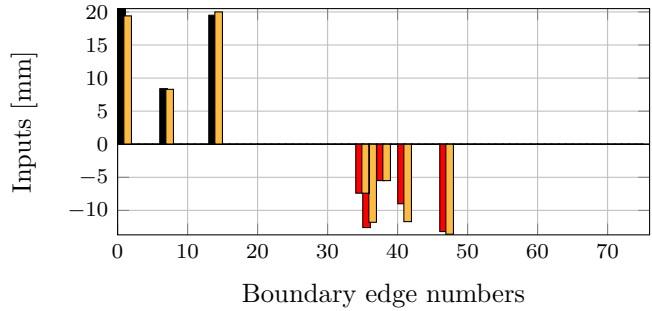
\begin{figure}[h!] 
\centering
%
%
\begin{tikzpicture}

\begin{axis}[%
width=7cm,
height=3cm,
at={(1.535433in,0.519685in)},
scale only axis,
area legend,
xmin=0,
xmax=76,
xlabel={Boundary edge numbers},
xmajorgrids,
xtick={0,10,20,30,40,50,60,70},
ytick={-10,-5,0,5,10,15,20},
ymin=-13.7,
ymax=20.5,
ylabel={Inputs [mm]},
ymajorgrids,
ticklabel style = {font=\footnotesize},
]
\addplot[ybar,bar width=0.1cm,bar shift=-0.04cm,draw=black,fill=black] plot table[row sep=crcr] {%
1	20.5\\
2	0\\
3	0\\
4	0\\
5	0\\
6	0\\
7	8.4\\
8	0\\
9	0\\
10	0\\
11	0\\
12	0\\
13	0\\
14	19.5\\
15	0\\
16	0\\
17	0\\
18	0\\
19	0\\
20	0\\
21	0\\
22	0\\
23	0\\
24	0\\
25	0\\
26	0\\
27	0\\
28	0\\
29	0\\
30	0\\
31	0\\
32	0\\
33	0\\
34	0\\
35	0\\
36	0\\
37	0\\
38	0\\
39	0\\
40	0\\
41	0\\
42	0\\
43	0\\
44	0\\
45	0\\
46	0\\
47	0\\
48	0\\
49	0\\
50	0\\
51	0\\
52	0\\
53	0\\
54	0\\
55	0\\
56	0\\
57	0\\
58	0\\
59	0\\
60	0\\
61	0\\
62	0\\
63	0\\
64	0\\
65	0\\
66	0\\
67	0\\
68	0\\
69	0\\
70	0\\
71	0\\
72	0\\
73	0\\
74	0\\
75	0\\
}; \label{refpos}
\addplot[ybar,bar width=0.1cm,bar shift=-0.04cm,draw=black,fill=red, opacity=1] plot table[row sep=crcr] {%
1	0\\
2	0\\
3	0\\
4	0\\
5	0\\
6	0\\
7	0\\
8	0\\
9	0\\
10	0\\
11	0\\
12	0\\
13	0\\
14	0\\
15	0\\
16	0\\
17	0\\
18	0\\
19	0\\
20	0\\
21	0\\
22	0\\
23	0\\
24	0\\
25	0\\
26	0\\
27	0\\
28	0\\
29	0\\
30	0\\
31	0\\
32	0\\
33	0\\
34	0\\
35	-7.4\\
36	-12.6\\
37	0\\
38	-5.5\\
39	0\\
40	0\\
41	-9\\
42	0\\
43	0\\
44	0\\
45	0\\
46	0\\
47	-13.2\\
48	0\\
49	0\\
50	0\\
51	0\\
52	0\\
53	0\\
54	0\\
55	0\\
56	0\\
57	0\\
58	0\\
59	0\\
60	0\\
61	0\\
62	0\\
63	0\\
64	0\\
65	0\\
66	0\\
67	0\\
68	0\\
69	0\\
70	0\\
71	0\\
72	0\\
73	0\\
74	0\\
75	0\\
}; \label{refneg}
\addplot[ybar,bar width=0.1cm,bar shift=0.04cm,draw=black,fill=Dandelion] 
plot table[row sep=crcr] {%
1	19.4\\
2	0\\
3	0\\
4	0\\
5	0\\
6	0\\
7	8.3\\
8	0\\
9	0\\
10	0\\
11	0\\
12	0\\
13	0\\
14	20\\
15	0\\
16	0\\
17	0\\
18	0\\
19	0\\
20	0\\
21	0\\
22	0\\
23	0\\
24	0\\
25	0\\
26	0\\
27	0\\
28	0\\
29	0\\
30	0\\
31	0\\
32	0\\
33	0\\
34	0\\
35	0\\
36	0\\
37	0\\
38	0\\
39	0\\
40	0\\
41	0\\
42	0\\
43	0\\
44	0\\
45	0\\
46	0\\
47	0\\
48	0\\
49	0\\
50	0\\
51	0\\
52	0\\
53	0\\
54	0\\
55	0\\
56	0\\
57	0\\
58	0\\
59	0\\
60	0\\
61	0\\
62	0\\
63	0\\
64	0\\
65	0\\
66	0\\
67	0\\
68	0\\
69	0\\
70	0\\
71	0\\
72	0\\
73	0\\
74	0\\
75	0\\
}; 
\addplot[ybar,bar width=0.1cm,bar shift=0.04cm,draw=black,fill=Dandelion] 
plot table[row sep=crcr] {%
1	0\\
2	0\\
3	0\\
4	0\\
5	0\\
6	0\\
7	0\\
8	0\\
9	0\\
10	0\\
11	0\\
12	0\\
13	0\\
14	0\\
15	0\\
16	0\\
17	0\\
18	0\\
19	0\\
20	0\\
21	0\\
22	0\\
23	0\\
24	0\\
25	0\\
26	0\\
27	0\\
28	0\\
29	0\\
30	0\\
31	0\\
32	0\\
33	0\\
34	0\\
35	-7.4\\
36	-11.8\\
37	0\\
38	-5.5\\
39	0\\
40	0\\
41	-11.7\\
42	0\\
43	0\\
44	0\\
45	0\\
46	0\\
47	-13.6\\\\-13.7\\
48	0\\
49	0\\
50	0\\
51	0\\
52	0\\
53	0\\
54	0\\
55	0\\
56	0\\
57	0\\
58	0\\
59	0\\
60	0\\
61	0\\
62	0\\
63	0\\
64	0\\
65	0\\
66	0\\
67	0\\
68	0\\
69	0\\
70	0\\
71	0\\
72	0\\
73	0\\
74	0\\
75	0\\
}; \label{sparse}
\end{axis} 
\end{tikzpicture}%
\caption{\label{fig:bar_ref_sparse} 
\ref{refpos} Reference inputs for shortening and 
\ref{refneg} for lengthening to steer the initial perturbed configuration  $r^{\mathrm{(ini)}}$ to the desired target configuration $r^{\mathrm{(des)}}$,  
 \ref{sparse} Computed sparse control inputs by Algorithm~\ref{alg:reweight}. }
\end{figure}

For comparison, both the fully actuated and the sparse input vectors are computed and their performance is compared in terms of the error norms $\|r^{(\mathrm{des})} - r^{(\mathrm{con})}\|_{Q_r}^2$ and $\|r^{(\mathrm{des})} - \rl^{(\mathrm{con})}\|_{Q_r}^2$, respectively. 
Note that in this case the solution is known to be sparse because of the definition of the initial condition and the target. 
Figure~\ref{fig:bar_ref_dense} shows the fully actuated control input vector $u$ computed by Algorithm~\ref{alg:SQP} and Figure~\ref{fig:bar_ref_sparse} shows the sparse input vector $\ul$ resulting from Algorithm~\ref{alg:reweight} with parameters $\tau = 10^{-4}$, $\epsilon = 10^{-8}$ and $\gamma = 0.3$. 
In the experiment, only the sparse input vector is applied to the prototype leading from the initial perturbed (blue) configuration, $r^{\mathrm{(ini)}}$, to the controlled (green) one, $\rl^{(\mathrm{con})}$. Because of the very good control performance, only one control iteration is done on the prototype. 
The fully actuated control input vector $u$ is not experimentally applied to the prototype system. However, in simulation, both $u$ and the sparse $\ul$ can be compared. With $r^{(\mathrm{con})}$ and $\rl^{(\mathrm{con})}$ being the minimizers of Problem $\mathcal{P}_{\mathrm{minE}}$ for the fully actuated $u$ and the sparse $\ul$, the error norms are $\|r^{(\mathrm{des})} - r^{(\mathrm{con})}\|_{Q_r}^2 = 1.669 \times 10^{-3}$ and $\|r^{(\mathrm{des})} - \rl^{(\mathrm{con})}\|_{Q_r}^2 = 1.883 \times 10^{-3}$, respectively. 

The experiments on the prototype show the following control performance. The measured data reveal that the error norm is decreased by $98.8 \%$, from ${\| r^{(\mathrm{des})} - r^{(\mathrm{ini})} \|^2_{Q_r}} = 1.55 \times 10^{-2}$ to ${\| r^{(\mathrm{des})} - \rl^{(\mathrm{con})} \|^2_{Q_r}} = 1.82 \times 10^{-4}$. 
The unweighted error norm of the deviations is decreased by $98.7 \%$ from ${\| r^{(\mathrm{des})} - r^{(\mathrm{ini})} \|^2_{2} = 2.21 \times 10^{-2}}$ to ${\| r^{(\mathrm{des})} - \rl^{(\mathrm{con})} \|^2_{2}} = 2.81 \times 10^{-4}$. 
The RMS-error, defined by $\| r^{(\mathrm{des})} - \rl^{(\mathrm{con})} \|_2 / n$, decreases by $88.5 \%$ from $\SI{0.134}{\mm}$ to $\SI{0.0154}{\mm}$. 

\begin{figure}
\vspace{0.0cm}
\includegraphics[width=0.6\columnwidth]{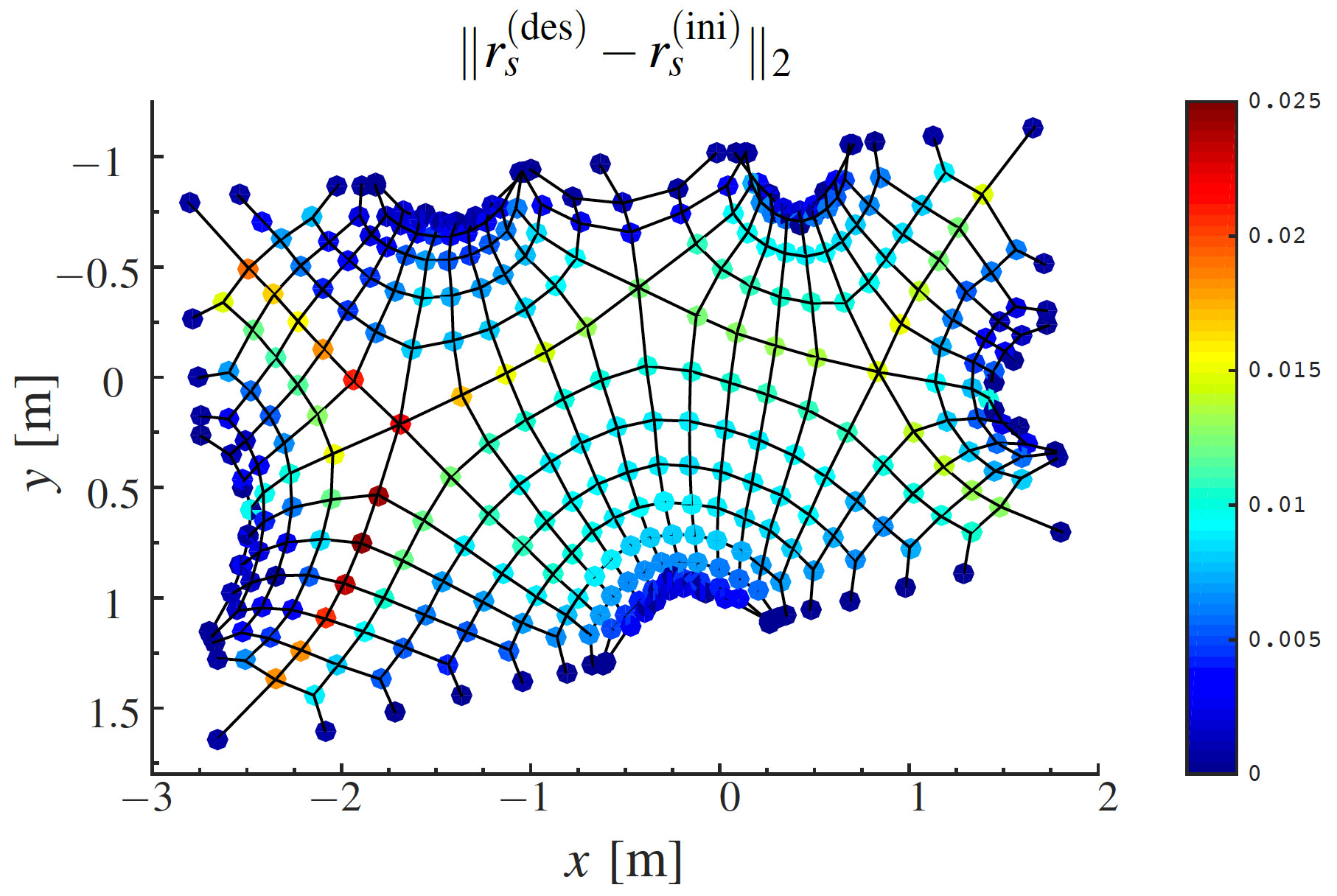} 
\caption{Control scenario with sparse actuation. Initial errors $\| r_s^{(\mathrm{des})} - r_s^{(\mathrm{ini})} \|_2$, the distance between the initial configuration $r_s^{(\mathrm{ini})}$ and the desired coordinates $r_s^{(\mathrm{des})}$, for each node $s = 1,...,n$. }
\label{fig:heat_rc1_bef}
\end{figure}

\begin{figure}
\vspace{0.0cm}
\includegraphics[width=0.6\columnwidth]{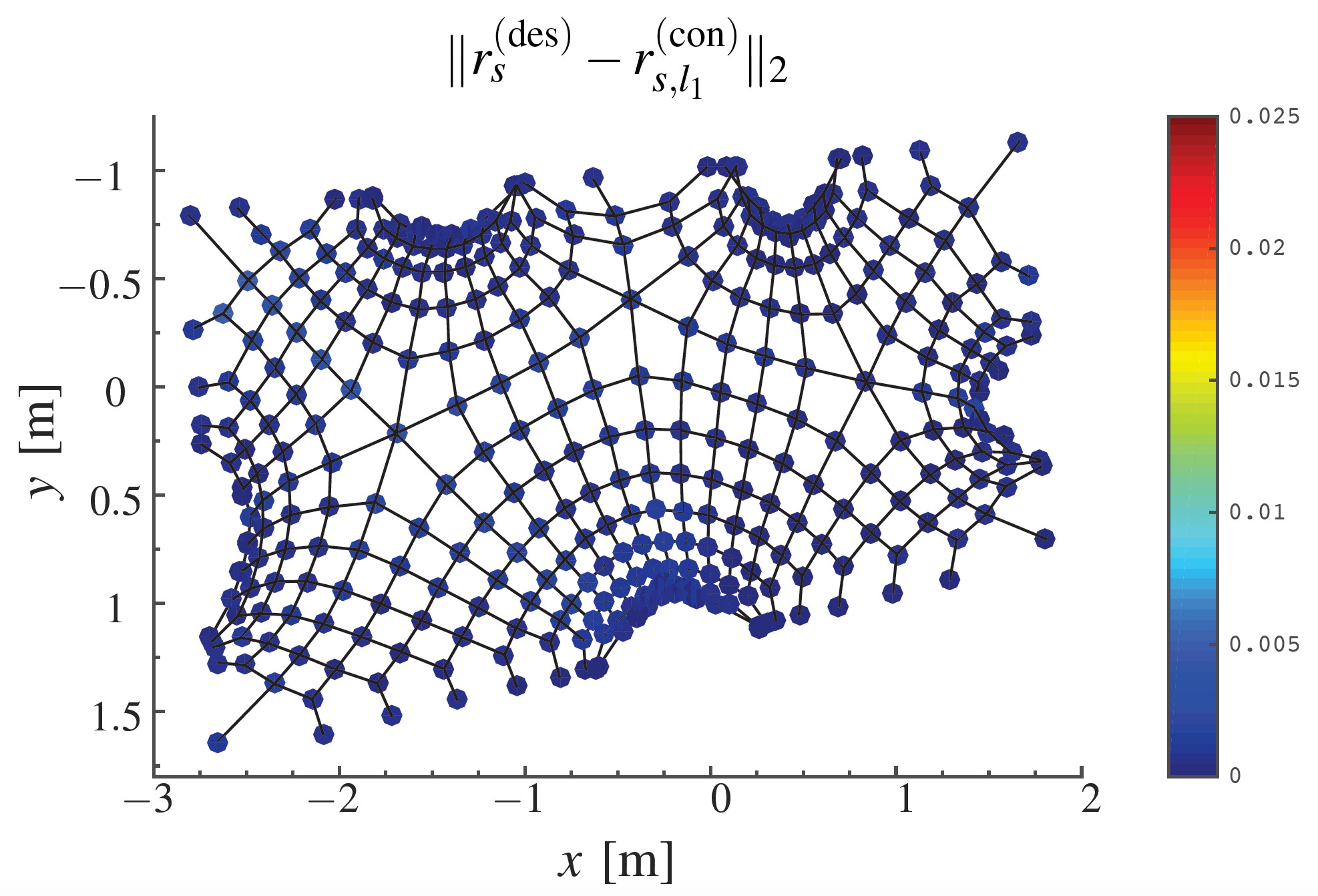}  \\
\caption{Control scenario with sparse actuation. Controlled errors $\| r_s^{(\mathrm{des})} - \rsl^{(\mathrm{con})} \|_2$, the distance between the controlled coordinates,  $\rsl^{(\mathrm{con})}$, and the desired coordinates, $r_s^{(\mathrm{des})}$. 
}
\label{fig:heat_rc1_aft}
\end{figure}

Figure~\ref{fig:heat_rc1_bef} shows the spatial distribution of the initial nodal position errors over the net. 
The spatial distribution of the controlled nodal position errors is shown in Figure~\ref{fig:heat_rc1_aft}. 

\begin{figure}
\centering
%
%
\definecolor{mycolor1}{rgb}{0.00000,0.44700,0.74100}%
\definecolor{mycolor2}{rgb}{0.85000,0.32500,0.09800}%
\begin{tikzpicture}
\begin{axis}[%
width=7.6cm, 
height=3cm, 
at={(1.535433in,0.519685in)},
scale only axis,
area legend,
xlabel={Deviations [mm]},
ymin=0,
ymax=160,
xmin=0,
xmax=28,
ytick={0,40,80,120,160},
xtick={ 0, 2, 4, 6, 8, 10, 12, 14, 16, 18, 20, 22, 24, 26, 28},
ticklabel style = {font=\footnotesize},
every axis y label/.style={at={(current axis.north west)},left=10mm, below=15mm, rotate=90}, 
ylabel={Count},
]
\addplot[ybar,bar width=0.08cm,bar shift=-0.02cm,draw=black,fill=Peach,opacity=0.5] plot table[row sep=crcr] {%
0	39\\
0.3	34\\
0.6	5\\
0.9	6\\
1.2	5\\
1.5	3\\
1.8	6\\
2.1	8\\
2.4	9\\
2.7	6\\
3.0	8\\
3.3	4\\
3.6	8\\
3.9	3\\
4.2	5\\
4.5	8\\
4.8	2\\
5.1	12\\
5.4	3\\
5.7	9\\
6.0	12\\
6.3	8\\
6.6	10\\
6.9	11\\
7.2	3\\
7.5	5\\
7.8	11\\
8.1	7\\
8.4	14\\
8.7	20\\
9	10\\
9.3	6\\
9.6	5\\
9.9	8\\
10.2	6\\
10.5	7\\
10.8	4\\
11.1	1\\
11.4	0\\
11.7	2\\
12	5\\
12.3	2\\
12.6	3\\
12.9	3\\
13.2	0\\
13.5	3\\
13.8	1\\
14.1	1\\
14.4	0\\
14.7	4\\
15.0	2\\
15.3	1\\
15.6	0\\
15.9	0\\
16.2	0\\
16.5	1\\
16.8	0\\
17.1	1\\
17.4	0\\
17.7	1\\
18	1\\
18.3	1\\
18.6	1\\
18.9	0\\
19.2	0\\
19.5	0\\
19.8	0\\
20.1	0\\
20.4	1\\
20.7	1\\
21	0\\
21.3	0\\
21.6	0\\
21.9	1\\
22.2	0\\
22.5	0\\
22.8	0\\
23.1	1\\
23.4	0\\
23.7	1\\
24	1\\
24.3	0\\
24.6	0\\
}; \label{before}
\addplot[ybar,bar width=0.08cm,bar shift=-0.02cm,draw=black,fill=mycolor1,opacity=0.5] plot table[row sep=crcr] {%
0	76\\
0.3	127\\
0.6	85\\
0.9	44\\
1.2	17\\
1.5	5\\
1.8	8\\
2.1	1\\
2.4	1\\
2.7	1\\
3	2\\
3.3	0\\
3.6	1\\
3.9	0\\
4.2	0\\
4.5	2\\
4.8	0\\
5.1	0\\
5.4	0\\
5.7	0\\
6	0\\
6.3	0\\
6.6	0\\
6.9	0\\
7.2	0\\
7.5	0\\
7.8	0\\
8.1	0\\
8.4	0\\
8.7	0\\
9	0\\
9.3	0\\
9.6	0\\
9.9	0\\
10.2	0\\
10.5	0\\
10.8	0\\
11.1	0\\
11.4	0\\
11.7	0\\
12	0\\
12.3	0\\
12.6	0\\
12.9	0\\
13.2	0\\
13.5	0\\
13.8	0\\
14.1	0\\
14.4	0\\
14.7	0\\
15	0\\
15.3	0\\
15.6	0\\
15.9	0\\
16.2	0\\
16.5	0\\
16.8	0\\
17.1	0\\
17.4	0\\
17.7	0\\
18	0\\
18.3	0\\
18.6	0\\
18.9	0\\
19.2	0\\
19.5	0\\
19.8	0\\
20.1	0\\
20.4	0\\
20.7	0\\
21	0\\
21.3	0\\
21.6	0\\
21.9	0\\
22.2	0\\
22.5	0\\
22.8	0\\
23.1	0\\
23.4	0\\
23.7	0\\
24	0\\
24.3	0\\
24.6	0\\
}; \label{after}
\end{axis}
\end{tikzpicture}%
\caption{Histograms of the distances  \ref{before} from the initial coordinates to the desired ones:  ${\| r_s^{(\mathrm{des})} - r_s^{(\mathrm{ini})}\|_2}$  \ref{after} from the controlled coordinates to the desired ones: ${\| r_s^{(\mathrm{des})} - \rsl^{(\mathrm{con})}\|_2}$, for all nodes $s = 1,...,n$.}
\label{fig:hist_r_c2_bef_aft}
\end{figure}
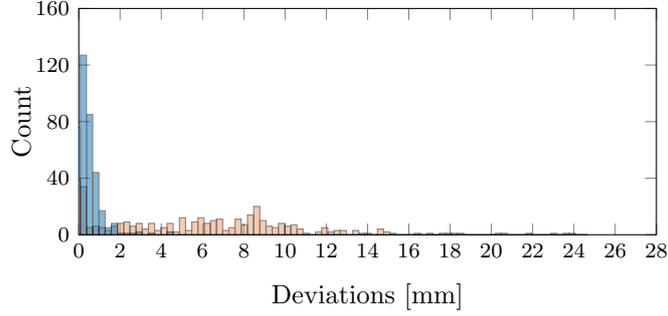

Figure~\ref{fig:hist_r_c2_bef_aft} shows the histogram of the statistical distribution of the measured deviation before and after the control on the prototype in terms of the Euclidean distances. The highest deviations can be seen in the $z$-coordinates, which are corrected from initial errors of more than $\SI{15}{\mm}$ to final errors of approximately $\SI{2}{\mm}$.

\section{Conclusions}
\label{sec:Conclusion}
A new control application has been presented for the control of the form of an innovative flexible formwork for concrete construction. This enables the precise construction of thin light-weight shell structures. The main component of the formwork, a pre-stressed cable net structure, needs to be controlled in order to ensure that the concreted shell precisely achieves the designed and optimized form and thus obtains its designed structural properties. 
Actuation of the system is possible by changing the lengths of the boundary edges, which are attached to a supporting rigid frame. 
The control algorithm minimizes the error norm between the built configuration of the net and the desired one. It is based on SQP with a guaranteed feasible step at each iteration. For practical application, sparse input vectors can be computed. 
A 1:4 prototype of the HiLo roof is used for the experimental assessment of the method. The experiments on the prototype system show very good control performance. The RMS-errors between desired and measured configurations before and after the control are decreased by $88.51 \%$ from $\SI{0.134}{\mm}$ to $\SI{0.0154}{\mm}$.   
For the construction of the HiLo roof on the NEST building, hardware improvements are planned in order to reduce the estimation errors of the nodal positions. 
As the roof has a larger span, the ratio between the number of interior nodal positions to be controlled and boundary edges available for control inputs will increase. Investigating a measure of controllability as a function of this ratio will be relevant. 
Furthermore, the pre-stressing forces will be higher and the material of the edges and of the supporting frame will be steel instead of the plastic used in the experimental investigation. Depending on the structural design, constraints on tension bounds in the edges might need to be introduced in order to prevent the material from being damaged.  

\section*{Acknowledgment}
The authors would like to thank Prof.\ Philippe Block, Dr.\ Andrew Liew and Dr.\ Tom Van Mele, from the Block Research Group, ETH Zurich, for the collaboration and the supervision of the experimental setup. Furthermore, many thanks to Jean-Marc Stadelmann for the actual construction of the prototype, and Dr.\ S{\'e}bastien Guillaume, from the Institute of Geodesy and Photogrammetry, ETH Zurich, for the measurements, Michael Lyrenmann from the NCCR Digital Fabrication, ETH Zurich, for taking the photographs, and Damian Frick from the Automatic Control Laboratory, ETH Zurich, for valuable technical discussions. 

\appendix

\subsection*{Proof of Lemma~\ref{lem:equalGNdirection}}
The following Corollary~\ref{cor:GN} is needed to prove Lemma~\ref{lem:equalGNdirection}. 
\begin{mycor} \label{cor:GN}
The Jacobian $\nabla_u \frF(u^\iter)$ in each iteration has singular values uniformly bounded away from zero, i.e., 
\begin{equation} \label{eq:Jacobian_assumption}
\exists ~ \nu > 0 ~\text{such that}~ \| \nabla_u \frF(u^\iter) \, \tilde{u}\| \geq \nu \,\|\tilde{u}\|\,, ~\forall \tilde{u} \in \mathbb{R}^{m_B},
\end{equation}
for all $u^\iter$ in a neighborhood of the bounded level set ${\mathcal{L} = \{u^\iter \vert \focpu(u^\iter) \leq \focpu(u^{0})\}}$, with $u^{0}$ being the starting point of the iteration.  
\end{mycor}
\begin{IEEEproof}
According to Theorem~\ref{the:implicitfunction}, we have 
${\nabla_u \frF(u) = - [\nabla_{r_F} h(r_F,u) ]^{-1} \nabla_u h(r_F,u)}$. 
Herein, the partial Jacobian ${\nabla_{r_F} h(r_F,u)}$ is invertible, because of Proposition~\ref{prop:Jacinv} and $\nabla_u h(r_F,u)$ has its singular values bounded away from zero because of Proposition~\ref{prop:Jacrank}. 
Moreover, both terms ${\nabla_{r_F} h(r_F,u)}$ and ${\nabla_u h(r_F,u)}$ have singular values bounded away from zero and bounded above, because first-order changes in both the nodal positions $r_F$ and the inputs $u$ have linearly independent non-zero, but bounded effects on the resulting forces at the free nodes, which completes the proof.
\end{IEEEproof}

\setcounter{mylem}{1} 
\begin{mylem} 
In each iteration $\iter$, the GN search direction $\Delta \uPocpu^\iter$ for $\Pocpu$ is equal to the partial GN search direction $\Delta u^\iter$ of $\PSQPGN$ in Algorithm~\ref{alg:SQP}. 
\end{mylem}
\begin{IEEEproof}
First, we reformulate the cost function $\fSQPGN(\Delta r_F^\iter)$ in \eqref{eq:GN_f} into 
\begin{equation} \label{eq:focpuGN}
\begin{aligned}
\focpuGN(\Delta \uPocpu^\iter) = \frac{1}{2} \left\Vert \frF(u^\iter)- {r}_{F}^{\mathrm{des}} + \nabla_u \frF(u^\iter) \Delta \uPocpu^\iter \right\Vert_{Q_r}^2. 
\end{aligned}
\end{equation}
We claim that at each feasible iterate, $\focpuGN(\Delta u^\iter)$ in \eqref{eq:focpuGN} is equal to $\fSQPGN(\Delta r_F^\iter)$ in \eqref{eq:GN_f}.  
First, we see that for any feasible iterate, it holds that $\frF(u^\iter) = r_F^\iter$. 
Moreover, with ${\nabla_u h(r_F^\iter,u^\iter) \, \Delta u^\iter + \nabla_{r_F} h(r_F^\iter,u^\iter) \, \Delta r_F^\iter = 0}$, we have ${\Delta r_F^\iter = - \left( \nabla_{r_F} h(r_F^\iter,u^\iter) \right)^{-1} \, \nabla_u h(r_F^\iter,u^\iter) \, \Delta u^\iter}$. With the Implicit Function Theorem, it holds that ${\Delta r_F^\iter = \nabla_u \frF(u^\iter) \, \Delta u^\iter}$ and thus $\focpuGN(\Delta u^\iter)$ and $\fSQPGN(\Delta r_F^\iter)$ are equal at the current point. 

As the constraints of $\PSQPGN$ satisfy the LICQ, 
the tangent cone of the nonlinear constraints $h(r_F,u)=0$ and the set of feasible linearized directions of $\PSQPGN$ are the same at the current point. 

Because of Lemma~\ref{lem:SQP}, the solution of $\PSQPGN$ is unique, 
and with the rank argument in Corollary~\ref{cor:GN}, it is easy to see from \eqref{eq:GN_searchdirection} that also $\Delta \uPocpu^\iter$ is unique. 
Therefore, the minimizer $\Delta \uPocpu^\iter$ for $\Pocpu$ is equal to the partial minimizer $\Delta u^\iter$ of $\PSQPGN$.  
\end{IEEEproof}

\addtolength{\textheight}{-0cm}
\bibliography{cablenetTCST}{}
\bibliographystyle{IEEEtran}

\end{document}